\documentclass[12pt]{amsart}
\usepackage{amsmath, amssymb, amsfonts}
\usepackage[all]{xypic}
\usepackage[pdftex]{graphicx}
\usepackage{upgreek}
\usepackage{mathtools}
\usepackage{stmaryrd}
\usepackage{bbm}

\theoremstyle{plain}
\newtheorem{theorem}{Theorem}[section]
\newtheorem{lemma}[theorem]{Lemma}
\newtheorem{corollary}[theorem]{Corollary}
\newtheorem{prop}[theorem]{Proposition}
\newtheorem{prop-def}[theorem]{Proposition / Definition}

\newtheorem{conj}[theorem]{Conjecture}

\theoremstyle{remark}

\newtheorem{remark}[theorem]{Remark}

\newtheorem{example}[theorem]{Example}

\newtheorem*{note*}{Note}
\newtheorem*{remark*}{Remark}
\newtheorem*{example*}{Example}

\theoremstyle{definition}
\newtheorem*{definition*}{Definition}
\newtheorem*{hypothesis*}{Hypothesis}
\newtheorem*{assumptions*}{Assumptions}
\newtheorem{definition}[theorem]{Definition}

\newcommand{\Z}{\mathbb{Z}}
\newcommand{\R}{\mathbb{R}}
\newcommand{\Q}{\mathbb{Q}}
\newcommand{\C}{\mathbb{C}}

\newcommand{\F}{\mathbb{F}}
\newcommand{\Aff}{\mathrm{Aff}}
\newcommand{\Aut}{\mathrm{Aut}}
\newcommand{\Gal}{\mathrm{Gal}}
\newcommand{\Tr}{\mathrm{Tr}}
\newcommand{\GL}{\mathrm{GL}}
\newcommand{\rank}{\mathrm{rank}}
\newcommand{\Nrd}{\mathrm{Nrd}}
\newcommand{\coker}{\mathrm{coker}}
\newcommand{\End}{\mathrm{End}}
\newcommand{\Hom}{\mathrm{Hom}}

\newcommand{\tors}{\mathrm{tors}}
\newcommand{\Irr}{\mathrm{Irr}}
\newcommand{\ind}{\mathrm{ind}}
\newcommand{\PMod}{\mathrm{PMod}}
\newcommand{\Spec}{\mathrm{Spec}}
\newcommand{\ram}{\mathrm{ram}}
\newcommand{\id}{\mathrm{id}}
\newcommand{\Sch}{\mathrm{Sch}}
\newcommand{\Char}{\mathrm{Char}}
\newcommand{\half}{{\textstyle \frac{1}{2}}}

\numberwithin{equation}{section}
\usepackage{tabto}

\newcommand{\et}{\mathrm{\acute{e}t}}
\newcommand{\syn}{\mathrm{syn}}
\newcommand{\perf}{\mathrm{perf}}
\newcommand{\infl}{\mathrm{infl}}

\newcommand{\Ext}{\mathrm{Ext}}
\newcommand{\reg}{\mathrm{reg}}
\newcommand{\cyc}{\mathrm{cyc}}
\newcommand{\nl}{\mathrm{nl}}

\DeclareFontFamily{U}{wncy}{}
    \DeclareFontShape{U}{wncy}{m}{n}{<->wncyr10}{}
    \DeclareSymbolFont{mcy}{U}{wncy}{m}{n}
    \DeclareMathSymbol{\Sha}{\mathord}{mcy}{"58}
    
\title[{On the $p$-adic Beilinson conjecture and the ETNC}]{On the $p$-adic Beilinson conjecture\\ and the equivariant Tamagawa number conjecture}

\author{Andreas Nickel}
\address{
Universit\"at Duisburg-Essen\\
Fakult\"at f\"ur Mathematik\\
Thea-Leymann-Stra\ss e 9\\
D-45127 Essen\\
Germany}
\email{andreas.nickel@uni-due.de}
\urladdr{https://www.uni-due.de/$\sim$hm0251/index.html}

\date{Version of 2nd October 2021}
\begin{document}

\begin{abstract}
Let $E/K$ be a finite Galois extension of totally real number fields
with Galois group $G$. Let $p$ be an odd prime and let $r>1$ be 
an odd integer. The $p$-adic Beilinson conjecture relates the values at $s=r$ of
$p$-adic Artin $L$-functions attached to the irreducible characters of $G$ to
those of corresponding complex Artin $L$-functions.
We show that this conjecture, the equivariant Iwasawa main conjecture
and a conjecture of Schneider imply the `$p$-part' of
the equivariant Tamagawa number conjecture
for the pair $(h^0(\Spec(E))(r), \Z[G])$. If $r>1$ is even we obtain a
similar result for Galois CM-extensions after restriction to `minus parts'.
\end{abstract}

\maketitle

\section{Introduction}\label{sec:introduction}
Let $E/K$ be a finite Galois extension of number fields with Galois group $G$
and let $r$ be an integer. The equivariant Tamagawa number conjecture (ETNC)
for the pair $(h^0(\Spec(E))(r), \Z[G])$ as formulated by
Burns and Flach \cite{MR1884523} asserts that a certain canonical element
$T\Omega(E/K,r)$ in the relative algebraic $K$-group $K_0(\Z[G], \R)$ vanishes.
This element relates the leading terms at $s=r$ of Artin $L$-functions to
natural arithmetic invariants. 

If $r=0$ this might be seen as a vast generalization
of the analytic class number formula for number fields, and refines
Stark's conjecture for $E/K$ as discussed by Tate in \cite{MR782485}
and the `Strong Stark conjecture' of Chinburg \cite[Conjecture 2.2]{MR724009}.
It is known to imply a whole bunch of conjectures such as 
Chinburg's `$\Omega_3$-conjecture' \cite{MR724009, MR786352}, 
the Rubin--Stark conjecture \cite{MR1385509}, Brumer's conjecture,
the Brumer--Stark conjecture (see \cite[Chapitre IV, \S 6]{MR782485}) 
and generalizations thereof due to Burns
\cite{MR2845620} and the author \cite{MR2976321}.
If $r$ is a negative integer, the ETNC refines a conjecture of Gross \cite{MR2154331}
and implies (generalizations of) the Coates--Sinnott
conjecture \cite{MR0369322} and a conjecture of Snaith \cite{MR2209286}
on annihilators of the higher $K$-theory of rings of integers (see \cite{MR2801311}).
If $r>1$ the ETNC likewise predicts constraints on the Galois module structure
of $p$-adic wild kernels \cite{wild-kernels}.

The functional equation of Artin $L$-functions suggests that the ETNC at $r$
and $1-r$ are equivalent. This is not known in general, 
but leads to a further conjecture which is sometimes
referred to as the local ETNC. Except for the validity of the local ETNC it
therefore suffices to consider the (global) ETNC for either odd or even integers $r$.
Note that the local ETNC is widely believed to be easier to settle. For instance,
the `global epsilon constant conjecture' of Bley and Burns \cite{MR2005875} 
measures the compatibility of the closely related `leading term conjectures' at
$s=0$ \cite{MR1863302} and $s=1$ \cite{MR2371375}
and is known to hold for arbitrary tamely ramified extensions
\cite[Corollary 7.7]{MR2005875} and also for certain weakly ramified 
extensions \cite{MR3520002}.

Now suppose that $E/K$ is a Galois extension of totally real number fields
and let $p$ be an odd prime.
If $r<0$ is odd Burns \cite{MR3294653} and the author \cite{MR3072281}
independently have shown that the `$p$-part' of  
the ETNC for the pair $(h^0(\Spec(E))(r), \Z[G])$ holds
provided that a certain Iwasawa $\mu_p$-invariant vanishes (which conjecturally
is always true). The latter condition is mainly present because
the equivariant Iwasawa main conjecture (EIMC) for totally real fields then holds
by independent work of Ritter and Weiss \cite{MR2813337} and of Kakde \cite{MR3091976}.

The case $r \geq 0$ is more subtle.
Burns and Venjakob \cite{MR2290587, MR2749572} (see also 
\cite[Corollary 2.8]{MR3294653}) proposed a strategy for proving the
$p$-part of the ETNC for the pair $(h^0(\Spec(E))(1), \Z[G])$. More precisely,
this special case of the ETNC is implied by the vanishing of the relevant
$\mu_p$-invariant,
Leopoldt's conjecture for $E$ at $p$ and the `$p$-adic Stark conjecture at $s=1$'.
The latter conjecture relates the leading terms at $s=1$ of the complex and $p$-adic
Artin $L$-functions attached to characters of $G$ by certain comparison periods.
Note that Burns and Venjakob actually assume these conjectures for all odd primes $p$
and then deduce the ETNC for the pair $(h^0(\Spec(E))(1), \Z[\half][G])$,
but their approach has recently been refined by Johnston and the author \cite{MR4111943}
so that one can indeed work prime-by-prime.

There are similar results on minus parts if $L/K$ is a Galois CM-extension with Galois
group $G$, i.e.~$K$ is totally real and $L$ is a totally complex quadratic extension
of a totally real field $L^+$.
Namely, if $r<0$ is even and $\mu_p$ vanishes, then the minus $p$-part
of the ETNC for the pair $(h^0(\Spec(L))(r), \Z[G])$ holds \cite{MR3294653, MR3072281}.
Burns \cite{MR4092926} recently proposed a strategy for proving the minus $p$-part
of the ETNC in the case $r=0$. In comparison with the strategy in the case $r=1$,
Leopoldt's conjecture is replaced with the conjectural non-vanishing of Gross's 
regulator \cite{MR656068}, and the $p$-adic Stark conjecture
is replaced with the `weak $p$-adic Gross--Stark conjecture' 
\cite[Conjecture 2.12b]{MR656068}
(now a theorem for linear characters by work of
Dasgupta, Kakde and Ventullo \cite{MR3866887}). 
For an approach that only relies upon 
the validity of the EIMC we refer the reader to \cite{MR2805422,MR3552493}.

The aim of this article is to propose a similar strategy in the remaining cases,
i.e.~we will consider the ETNC for Tate motives $h^0(\Spec(L))(r)$ where $r>1$
and $L$ is a CM-field.
Note that we can treat all integers $r>1$ simultaneously as the `plus $p$-part'
of the ETNC for the pair $(h^0(\Spec(L))(r), \Z[G])$ naturally identifies with
the corresponding conjecture for the extension $L^+/K$ of totally real fields.
We show that the $p$-adic Beilinson conjecture at $s=r$, 
a conjecture of Schneider \cite{MR544704}
and the EIMC imply the plus (resp.~minus)
$p$-part of the ETNC for the pair $(h^0(\Spec(L))(r), \Z[G])$ if $r$ is
odd (resp.~even). 

We follow the formulation of the $p$-adic Beilinson conjecture in \cite{MR2520465}.
It relates the values at $s=r$ of the complex and $p$-adic Artin $L$-functions
by certain comparison periods involving Besser's syntomic regulator \cite{MR1809626}.
For absolutely abelian extensions variants of the $p$-adic Beilinson conjecture 
have been formulated and proved by Coleman \cite{MR674400}, Gros \cite{MR1031903, MR1248079}
and Kolster and Nguyen Quang Do \cite{MR1632798}.
Thus the $p$-adic Beilinson conjecture holds for absolutely abelian characters
(see \S \ref{subsec:absolutely} for a precise statement).

Let us compare our approach to the earlier work mentioned above.
The formulation of both the $p$-adic Beilinson conjecture and the $p$-adic
Stark conjecture involves the choice of
a field isomorphism $j: \C \simeq \C_p$. We show in \S \ref{subsec:Gross-relation}
that the $p$-adic Beilinson conjecture does not depend upon this choice if and only if
a conjecture of Gross \cite{MR2154331} holds. The latter 
is revisited in \S \ref{subsec:Gross} and might be seen
as a higher analogue of Stark's conjecture; a similar result
in the case $r=1$ has recently been established by Johnston and the author
in \cite{MR4111943}. In both cases the independence of $j$ is therefore equivalent
to the rationality part of the appropriate special case of the ETNC. 
This eventually allows us to establish a
prime-by-prime descent result analogous to \cite[Theorem 8.1]{MR4111943}.

In a little more detail, we formulate 
conjectural `higher refined $p$-adic class number formulae'
analogous to \cite[Conjecture 3.5]{MR4092926} (where $r=0$),
and show that these follow from the EIMC and Schneider's conjecture
in \S \ref{subsec:higher-CNF}. 
Here, as will be shown in \S \ref{subsec:Schneider},
the latter conjecture ensures that the relevant complexes are semisimple
at all Artin characters as Leopoldt's conjecture does in the case $r=1$
and the non-vanishing of Gross's regulator does in the case $r=0$.
This is a necessary condition in order to apply the descent formalism
of Burns and Venjakob \cite{MR2749572}.
A second condition is the vanishing of the aforementioned Iwasawa $\mu_p$-invariant,
but given recent progress of Johnston and the author \cite{MR3749195, abelian-MC}
on the EIMC without assuming $\mu_p=0$,
we wish to circumvent this hypothesis. For this purpose, 
we develop a different descent argument that makes no use of this assumption,
but requires a more delicate analysis of the relevant complexes.
The higher refined $p$-adic class number formula at $s=r$ may then be combined 
with the $p$-adic Beilinson conjecture at $s=r$ to deduce the 
plus, respectively minus, $p$-part
of the ETNC for the pair $(h^0(\Spec(L))(r), \Z[G])$
in \S \ref{subsec:application}. For this, it is crucial to relate
Besser's syntomic regulators to Soul\'e's $p$-adic Chern class maps \cite{MR553999}
and the Bloch--Kato exponential maps \cite{MR1086888}
that appear in the formulation of the ETNC. This is carried out
in \S \ref{sec:comparison} (in particular see Proposition \ref{prop:crucial}).
The formulation of the ETNC that is most suitable for our purposes is a reformulation
due to the author \cite{wild-kernels}. This has primarily been introduced
in order to construct (conjectural) annihilators of $p$-adic wild kernels.

Our prime example are totally real Galois extensions $E/\Q$ with Galois
group isomorphic to $\Aff(q)$, where $q = \ell^n$ is a prime power
and $\Aff(q)$ denotes the group of affine transformations on the finite field
$\F_q$ with $q$ elements. We show that Gross's conjecture holds in 
this case (Theorem \ref{thm:Gross-cases} (iv)). Moreover, the relevant
cases of the EIMC hold unconditionally by recent work of Johnston and the author
\cite{MR3749195} (see also \cite{abelian-MC})
and the $p$-adic Beilinson conjecture reduces to the case of the trivial
extension $E^H/E^H$ where $H$ denotes the subgroup $\GL_1(\F_q)$ of $\Aff(q)$.
See Example \ref{ex:Aff(q)} for more details.

Finally, we note that the ETNC for the pair $(h^0(\Spec(L))(r), \Z[G])$
has been verified for any integer $r$ whenever $L$ is abelian over the rationals
by work of Burns, Greither and Flach \cite{MR1992015, MR2863902, MR2290586}.
However, if $r>1$ and $L$ is not absolutely abelian, then we are not aware of
any previous (conditional) results that establish the ($p$-part of the) ETNC
for the pair $(h^0(\Spec(L))(r), \Z[G])$.

\subsection*{Acknowledgements}
The author acknowledges financial support provided by the 
Deutsche Forschungsgemeinschaft (DFG) 
within the Heisenberg programme (project number 334383116).

\subsection*{Notation and conventions}
All rings are assumed to have an identity element and all modules are assumed
to be left modules unless otherwise stated. 
Unadorned tensor products will always denote
tensor products over $\Z$.
For a ring $\Lambda$ we write
$\zeta(\Lambda)$ for its center and $\Lambda^{\times}$ for the group of units in
$\Lambda$. 
For every field $F$ we fix a separable closure $F^c$ of $F$ and write
$G_F := \Gal(F^c/F)$ for its absolute Galois group.
If $n>0$ is an integer coprime to the characteristic of $F$, we let
$\zeta_n$ denote a primitive $n$-th root of unity in $F^c$.

A finite Galois extension of totally real number fields will usually be denoted
by $E/K$, whereas $L/K$ denotes an arbitrary Galois extension of number fields.
Galois CM-extensions will usually be denoted by $L/K$ as well.

\section{Algebraic Preliminaries}

\subsection{Derived categories and Galois cohomology} \label{sec:derived}
Let $\Lambda$ be a noetherian ring and let $\PMod(\Lambda)$ be the category of all finitely
generated projective $\Lambda$-modules. We write $\mathcal D(\Lambda)$ for the derived category
of $\Lambda$-modules and $\mathcal C^b(\PMod(\Lambda))$ for the category of bounded complexes
of finitely generated projective $\Lambda$-modules.
Recall that a complex of $\Lambda$-modules is called perfect if it is isomorphic in $\mathcal D(\Lambda)$
to an element of $\mathcal C^b(\PMod(\Lambda))$.
We denote the full triangulated subcategory of $\mathcal D(\Lambda)$
comprising perfect complexes by $\mathcal D^{\perf}(\Lambda)$.

If $M$ is a $\Lambda$-module and $n$ is an integer, we write $M[n]$ for the complex
\[
\cdots  \longrightarrow 0 \longrightarrow M \longrightarrow 0 \longrightarrow \cdots
\]
where $M$ is placed in degree $-n$. Note that this is compatible
with the usual shift operator on cochain complexes.

Let $L$ be an algebraic extension of the number field $K$.
For a finite set $S$ of places of $K$ containing the set $S_{\infty}$
of all archimedean places
we let $G_{L,S}$ be the Galois group over $L$ of the maximal extension of $L$
that is unramified outside $S(L)$; here we write $S(L)$
for the set of places of $L$ lying above those in $S$.
We let $\mathcal O_{L,S}$ be the ring of $S(L)$-integers in $L$.
For any topological $G_{L,S}$-module $M$
we write $R\Gamma(\mathcal O_{L,S}, M)$ for the complex
of continuous cochains of $G_{L,S}$ with coefficients in $M$.
If $F$ is a field and $M$ is a topological $G_{F}$-module, we likewise define
$R\Gamma(F,M)$ to be the complex
of continuous cochains of $G_{F}$ with coefficients in $M$.

If $F$ is a global or a local field of characteristic zero, and $M$ is a discrete or
a compact $G_F$-module, then for $r \in \Z$ we denote the $r$-th Tate twist of $M$ by $M(r)$.
Now fix a prime $p$ and suppose that $S$ also contains all $p$-adic places of $K$. Then
for each integer $i$ the cohomology group in degree $i$ of
$R\Gamma(\mathcal O_{L,S}, \Z_p(r))$ naturally identifies with
$H^i_{\et} (\mathcal{O}_{L,S}, \Z_p(r))$, the $i$-th \'etale cohomology group of the affine scheme
$\Spec(\mathcal{O}_{L,S})$ with coefficients in the \'etale $p$-adic sheaf $\Z_p(r)$.
We set $H^i_{\et} (\mathcal{O}_{L,S}, \Q_p(r)) := \Q_p \otimes_{\Z_{p}}
H^i_{\et} (\mathcal{O}_{L,S}, \Z_p(r))$.

\subsection{Representations and characters of finite groups}\label{subsec:reps-and-chars-of-finite-groups}
Let $G$ be a finite group and let $F$ be a field of characteristic zero. 
We write $R_{F}^{+}(G)$ for the set of characters attached to finite-dimensional $F$-valued representations of $G$,
and $R_{F}(G)$ for the ring of virtual characters generated by $R_{F}^{+}(G)$.  
Moreover, we denote the subset of irreducible characters in $R_{F}^{+}(G)$
and the ring of $F$-valued virtual characters of $G$
by $\Irr_{F}(G)$ and $\Char_{F}(G)$, respectively. 

For a subgroup $H$ of $G$ and $\psi \in R_{F}^{+}(H)$ we write $\ind_{H}^{G} \psi \in R_{F}^{+}(G)$ for the induced character;
for a normal subgroup $N$ of $G$ and $\chi \in R_{F}^{+}(G/N)$ we write $\infl_{G/N}^{G} \chi \in R_{F}^{+}(G)$ for the inflated character.
For $\sigma \in \Aut(F)$ and $\chi \in \Char_{F}(G)$ we set 
$\chi^{\sigma} := \sigma \circ \chi$ and note that this defines a group action
from the left even though we write exponents on the right of $\chi$.
We denote the trivial character of $G$ by $\mathbbm{1}_{G}$.

\subsection{$\chi$-twists}
Let $G$ be a finite group and let $F$ be a field of characteristic zero. 
If $M$ is a $\Z[G]$-module we let $M^G$ be the maximal submodule of $M$
upon which $G$ acts trivially. Likewise we write $M_G$ for the maximal
quotient module with trivial $G$-action. 
For any $\chi \in R_{F}^{+}(G)$ we fix a (left) $F[G]$-module $V_{\chi}$ 
with character $\chi$.
For any $F[G]$-module $M$ and any $\alpha \in \End_{F[G]}(M)$ we write
$M^{\chi}$ for the $F$-vector space 
\[
\Hom_{F[G]}(V_{\chi}, M) 
\simeq \Hom_{F}(V_{\chi}, M)^{G} 
\] 
and $\alpha^{\chi}$ for the induced map $(f \mapsto \alpha \circ f) \in \End_{F}(M^{\chi})$.
We note that $\det_{F}(\alpha^{\chi})$ is independent of the choice of $V_{\chi}$.
The following is \cite[Lemma 2.1]{MR4111943} and very similar to 
\cite[Chapitre I, 6.4]{MR782485}.

\begin{lemma}\label{lem:properties-of-rho-parts-functor}
Let $M$ be an $F[G]$-module and let $\alpha \in \End_{F[G]}(M)$.
Let $H$ be a subgroup of $G$ and let $M |_{H}$ denote $M$ considered as an $F[H]$-module.
Let $N$ be a normal subgroup of $G$ and view $M^N$ as an $F[G/N]$-module
in the obvious way.
\begin{enumerate}
\item If $\chi_{1},\chi_{2} \in R_{F}^{+}(G)$ then  
${\det}_{F}(\alpha^{\chi_{1} + \chi_{2}}) = {\det}_{F}(\alpha^{\chi_{1}})  {\det}_{F}(\alpha^{\chi_{2}})$.
\item If $\chi \in R^{+}_{F}(H)$ then $M^{\ind^{G}_{H} \chi} \simeq (M |_{H})^{\chi}$ and 
${\det}_{F}(\alpha^{\ind^{G}_{H} \chi}) = {\det}_{F}(\alpha^{\chi})$.
\item If  $\chi \in R^{+}_{F}(G/N)$ then 
$M^{\infl^{G}_{G/N} \chi} \simeq (M^{N})^{\chi}$ and 
${\det}_{F}(\alpha^{\infl^{G}_{G/N} \chi}) = {\det}_{F}( (\alpha |_{M^{N}})^{\chi})$.
\end{enumerate} 
\end{lemma}

Let $p$ be a prime. For each $\chi \in \Irr_{\C_{p}}(G)$ we fix a subfield
$F_{\chi}$ of $\C_p$ which is both Galois and of finite degree over $\Q_p$
and such that $\chi$ can be realized over $F_{\chi}$.
We write $e_{\chi} := |G|^{-1} \chi(1)
\sum_{g \in G} \chi(g) g^{-1}$ for the associated primitive central idempotent
in $\C_p[G]$ and choose an indecomposable idempotent $f_{\chi}$ of 
$F_{\chi}[G]e_{\chi}$. Let $\mathcal{O}_{\chi}$ be the ring of integers
in $F_{\chi}$ and choose a maximal $\mathcal{O}_{\chi}$-order
$\mathcal{M}_{\chi}$ in $F_{\chi}[G]$ containing $f_{\chi}$.
Then $T_{\chi} := f_{\chi} \mathcal{M}_{\chi}$ is an $\mathcal{O}_{\chi}$-free
right $\mathcal{O}_{\chi}[G]$-module.

For any (left) $\Z_p[G]$-module $M$ we define a (left) 
$\mathcal{O}_{\chi}[G]$-module $M[\chi] := T_{\chi} \otimes_{\Z_{p}} M$,
where $g \in G$ acts upon $t \otimes m$ by the rule $g(t \otimes m) = 
t g^{-1} \otimes gm$ for all $t \in T_{\chi}$ and $m \in M$.
We define $\mathcal{O}_{\chi}$-modules $M^{(\chi)} := M[\chi]^G$
and $M_{(\chi)} := M[\chi]_G \simeq T_{\chi} \otimes_{\Z_p[G]} M$.
We thereby obtain left, respectively right exact functors
$M \mapsto M^{(\chi)}$ and $M \mapsto M_{(\chi)}$ from the category of
$\Z_p[G]$-modules to the category of $\mathcal{O}_{\chi}$-modules.
Note that there is an isomorphism
$F_{\chi} \otimes_{\mathcal{O}_{\chi}} M^{(\chi)} \simeq 
(F_{\chi} \otimes_{\Z_{p}} M)^{\chi}$ for every finitely generated
$\Z_p[G]$-module $M$. 

Since multiplication by the trace $\Tr_G := \sum_{g \in G} g$ gives rise to an isomorphism
$P_{(\chi)} \simeq P^{(\chi)}$ for each projective $\Z_p[G]$-module $P$
(in fact for each cohomologically trivial $G$-module $P$), these functors
extend to naturally isomorphic exact functors $\mathcal{D}^{\perf}(\Z_p[G])
\rightarrow \mathcal{D}^{\perf}(\mathcal{O}_{\chi})$
(and $\mathcal{D}^{\perf}(\Q_p[G])
\rightarrow \mathcal{D}^{\perf}(F_{\chi})$).

\begin{lemma} \label{lem:chi-twists}
Let $\chi \in \Irr_{\C_{p}}(G)$ and let $a \leq b$ be integers.
If $C^{\bullet} \in \mathcal{D}^{\perf}(\Z_p[G])$ is acyclic outside $[a,b]$,
then $C^{\bullet}_{(\chi)}$ is also acyclic outside $[a,b]$ and there are 
natural isomorphisms of $\mathcal{O}_{\chi}$-modules
\[
H^a(C^{\bullet}_{(\chi)}) \simeq H^a(C^{\bullet})^{(\chi)}
\qquad \mbox{and} \qquad
H^b(C^{\bullet}_{(\chi)}) \simeq H^b(C^{\bullet})_{(\chi)}.
\]
For $C^{\bullet} \in \mathcal{D}^{\perf}(\Q_p[G])$ we have isomorphisms
$H^i(C^{\bullet}_{(\chi)}) \simeq H^i(C^{\bullet})_{(\chi)}
\simeq H^i(C^{\bullet})^{(\chi)}$ for every $i \in \Z$.
\end{lemma}

\begin{proof}
Since (finitely generated) $\Q_p[G]$-modules are cohomologically trivial,
the functors $M \mapsto M^{(\chi)}$ and $M \mapsto M_{(\chi)}$
are naturally isomorphic exact functors on the category of finitely generated
$\Q_p[G]$-modules. The final assertion of the lemma is therefore clear.

Now suppose that $C^{\bullet} \in \mathcal{D}^{\perf}(\Z_p[G])$ 
is acyclic outside $[a,b]$.
If $b-a \leq 1$ the claim is \cite[Lemma 5.1]{MR4092926}. We repeat
the short argument for convenience. Choose a complex $A \rightarrow B$
of cohomologically trivial $\Z_p[G]$-modules that is isomorphic to $C^{\bullet}$
in $\mathcal{D}(\Z_p[G])$. Here $A$ and $B$ are placed in degrees $a$ and $a+1$, 
respectively. Then we obtain a commutative diagram of $\mathcal{O}_{\chi}$-modules
\[ \xymatrix{
	& & A_{(\chi)} \ar[r] \ar[d]^{\simeq} & B_{(\chi)} \ar[d]^{\simeq} \ar[r] &
	H^{a+1}(C^{\bullet})_{(\chi)} \ar[r] & 0 \\
	0 \ar[r] & H^a(C^{\bullet})^{(\chi)} \ar[r] & A^{(\chi)} \ar[r] & 
	B^{(\chi)} & &
}\]
which implies the claim. If $b-a \geq 2$ 
we choose a complex $P^{\bullet} \in \mathcal{C}^b(\PMod(\Z_p[G]))$ 
that is isomorphic to $C^{\bullet}$
in $\mathcal{D}(\Z_p[G])$ and
consider the exact sequence of perfect complexes
\[ 
0 \longrightarrow \tau_{\geq b-1} P^{\bullet} \longrightarrow P^{\bullet}
\longrightarrow \tau_{\leq b-2} P^{\bullet} \longrightarrow 0,
\]
where $\tau_{\geq b-1}$ and $\tau_{\leq b-2}$ denote naive truncation.
Note that the complexes $\tau_{\geq b-1} P^{\bullet}$ and
$\tau_{\leq b-2} P^{\bullet}$ are acyclic outside $[b-1,b]$ and
$[a,b-2]$, respectively.
It follows by induction that $C^{\bullet}_{(\chi)}$ is acyclic outside $[a,b]$
and, since $H^b(C^{\bullet}) = H^b(\tau_{\geq b-1} P^{\bullet})$, 
that we have an isomorphism
$H^b(C^{\bullet}_{(\chi)}) \simeq H^b(C^{\bullet})_{(\chi)}$.
If $b-a \geq 3$ then we likewise have that
$H^a(C^{\bullet}) = H^a(\tau_{\leq b-2} P^{\bullet})$
and we may again conclude by induction that we have an isomorphism
$H^a(C^{\bullet}_{(\chi)}) \simeq H^a(C^{\bullet})^{(\chi)}$. If $b-a=2$ we may
alternatively consider the exact sequence of
perfect complexes
\[ 
0 \longrightarrow \tau_{\geq b} P^{\bullet} \longrightarrow P^{\bullet}
\longrightarrow \tau_{\leq b-1} P^{\bullet} \longrightarrow 0
\]
and deduce as above.
\end{proof}

\section{The $p$-adic Beilinson conjecture}

\subsection{Setup and notation} \label{sec:setup}
Let $L/K$ be a finite Galois extension of number fields with Galois group $G$.
For any place $v$ of $K$ we choose a place $w$ of $L$ above $v$ and write $G_w$ and $I_w$ for
the decomposition group and inertia subgroup of $L/K$ at $w$, respectively.
We denote the completions of $L$ and $K$ at $w$ and $v$ by $L_w$ and $K_v$, respectively,
and identify the Galois group of the extension $L_w / K_v$ with $G_w$.
For each non-archimedean place $w$
we let $\mathcal{O}_w$ be the ring of integers in $L_w$.
We identify $\overline{G_w} := G_w / I_w$ with the Galois group of the 
corresponding residue field
extension which we denote by $L(w) / K(v)$. Finally, we let $\phi_w \in \overline{G_w}$ be the Frobenius
automorphism, and we denote the cardinality of $K(v)$ by $N(v)$.
We let $S$ be a finite set of places of $K$ containing the set $S_{\infty}$
of archimedean places. If a prime $p$ is fixed, we will usually assume that
the set $S_p$ of all $p$-adic places is also contained in $S$.

By a Galois CM-extension of number fields we shall mean a finite Galois extension $L/K$
such that $K$ is totally real and $L$ is a CM-field. Thus complex conjugation
induces a unique automorphism $\tau$ in the center of $G$ and we denote
the maximal totally real subfield of $L$ by $L^+$. Then $L^+/K$ is also Galois
with group $G^+ := G / \langle \tau \rangle$.

\subsection{Higher $K$-theory}
For an integer $n \geq 0$  and a ring $R$ we write
$K_n(R)$ for the Quillen $K$-theory of $R$. In the cases $R = \mathcal O_{L,S}$ and $R=L$
the groups $K_n(\mathcal O_{L,S})$ and $K_n(L)$ are equipped with a natural $G$-action and for every
integer $r>1$ the inclusion $\mathcal O_{L,S} \subseteq L$ induces an isomorphism of $\Z[G]$-modules
\begin{equation} \label{eqn:odd-iso}
K_{2r-1}(\mathcal O_{L,S}) \simeq K_{2r-1}(L).
\end{equation}
Moreover, if $S'$ is a second finite set of places of $K$ containing $S$, then for every $r>1$ there is a natural exact sequence
of $\Z[G]$-modules
\begin{equation} \label{eqn:even-ses}
0 \longrightarrow K_{2r}(\mathcal O_{L,S}) \longrightarrow K_{2r}(\mathcal O_{L,S'}) \longrightarrow
\bigoplus_{w \in S'(L) \setminus S(L)} K_{2r-1}(L(w)) \longrightarrow 0.
\end{equation}
Both results \eqref{eqn:odd-iso}
and \eqref{eqn:even-ses} are due to Soul\'{e} \cite{MR553999};
see \cite[Chapter V, Theorem 6.8]{MR3076731}. We also note that sequence \eqref{eqn:even-ses}
remains left-exact in the case $r=1$.
The structure of the finite $\Z[\overline{G_w}]$-modules
$K_{2r-1}(L(w))$ has been determined by Quillen \cite{MR0315016} (see also
\cite[Chapter IV, Theorem 1.12 and Corollary 1.13]{MR3076731}) to be
\begin{equation} \label{eqn:K-theory_finitefield}
K_{2r-1}(L(w)) \simeq \Z[\overline{G_w}] / (\phi_w - N(v)^r).
\end{equation}

\subsection{The regulators of Borel and Beilinson} \label{sec:Borel}
Let $\Sigma(L)$ be the set of embeddings of $L$ into the complex numbers; we then have
$|\Sigma(L)| = r_1 + 2 r_2$, where $r_1$ and $r_2$ are the number of real embeddings and the number
of pairs of complex embeddings of $L$, respectively. 
For an integer $k \in \Z$ we define a finitely generated $\Z$-module
\[
H_k(L) := \bigoplus_{\Sigma(L)} (2 \pi i)^{-k} \Z
\]
which is endowed with a natural $\Gal(\C / \R)$-action, diagonally on $\Sigma(L)$ and on $(2\pi i)^{-k}$.
The invariants of $H_k(L)$ under this action will be denoted by $H_k^+(L)$, and it is easily seen that we have
\begin{equation} \label{eqn:rank-Betti}
d_k := \rank_{\Z}(H_{1-k}^+(L)) = \left\{ \begin{array}{lll}
r_1 + r_2 & \mbox{ if } & 2 \nmid k\\
r_2  & \mbox{ if } & 2 \mid k.
\end{array}\right.
\end{equation}
The action of $G$ on $\Sigma(L)$ endows $H_k^+(L)$ with a natural $G$-module structure.

Let $r>1$ be an integer. Borel \cite{MR0387496} has shown that the even $K$-groups
$K_{2r-2}(\mathcal{O}_L)$ (and thus $K_{2r-2}(\mathcal{O}_{L,S})$ for any $S$ as above
by \eqref{eqn:even-ses} and \eqref{eqn:K-theory_finitefield}) are finite,
and that the odd $K$-groups $K_{2r-1}(\mathcal{O}_L)$ are 
finitely generated abelian groups of
rank $d_r$. More precisely, for each $r>1$ Borel constructed an equivariant regulator map
\begin{equation} \label{eqn:regulator}
\rho_r^{\mathrm{Bor}}: K_{2r-1}(\mathcal{O}_L) \longrightarrow H_{1-r}^+(L) \otimes \R
\end{equation}
with finite kernel. Its image is a full lattice in $H_{1-r}^+(L) \otimes \R$.
The covolume of this lattice is called the Borel regulator and will be 
denoted by $R_r^{\mathrm{Bor}}(L)$.
Moreover, Borel showed that
\begin{equation} \label{eqn:Borels-rationality}
\frac{\zeta_L^{\ast}(1-r)}{R_r^{\mathrm{Bor}}(L)} \in \Q^{\times},
\end{equation}
where $\zeta_L^{\ast}(1-r)$ denotes the leading term at $s=1-r$ of 
the Dedekind zeta function $\zeta_L(s)$
attached to the number field $L$.

In the context of the ETNC, however, it is more natural to work with 
Beilinson's regulator map \cite{MR760999}.
By a result of Burgos Gil \cite{MR1869655} Borel's regulator map is twice the
regulator map of Beilinson. Hence we will work with $\rho_r := \half \rho_r^{\mathrm{Bor}}$
in the following.

\begin{remark} \label{rem:regulators}
We will sometimes refer to \cite{wild-kernels} where we have worked with
Borel's regulator map. However, if we are interested in rationality questions
or in verifying the $p$-part of the ETNC for an odd prime $p$, the factor $2$
essentially plays no role. In contrast, the $p$-adic Beilinson conjecture below
predicts an equality of two numbers in $\C_p$ so that this factor indeed
matters.
\end{remark}

\subsection{The Quillen--Lichtenbaum Conjecture}
Fix an odd prime $p$ and assume that $S$ contains $S_{\infty}$ and the set $S_p$ of
all $p$-adic places of $K$. Then for any integer $r > 1$ and $i=1,2$ Soul\'e \cite{MR553999}
has constructed canonical $G$-equivariant $p$-adic Chern class maps
\[
 ch^{(p)}_{r, i}: K_{2r-i}(\mathcal{O}_{L,S}) \otimes \Z_p \longrightarrow H^i_{\et} (\mathcal{O}_{L,S}, \Z_p(r)).
\]
We need the following deep result.

\begin{theorem}[Quillen--Lichtenbaum Conjecture] \label{thm:Quillen--Lichtenbaum}
 Let $p$ be an odd prime. Then for every integer $r > 1$ and $i=1,2$ the
 $p$-adic Chern class maps $ch^{(p)}_{r, i}$ are isomorphisms.
\end{theorem}

\begin{proof}
 Soul\'e \cite{MR553999} proved surjectivity. Building on work of Rost and Voevodsky, Weibel
 \cite{MR2529300} completed the proof of the Quillen--Lichtenbaum Conjecture.
\end{proof}

Let $p$ be a prime. For an integer $n \geq 0$ and a ring $R$ we write
$K_n(R; \Z_p)$ for the $K$-theory of $R$ with coefficients in $\Z_p$.
The following result is due to Hesselholt and Madsen \cite{MR1998478}.

\begin{theorem} \label{thm:local-Quillen--Lichtenbaum}
 Let $p$ be an odd prime and let $w$ be a finite place of $L$. Then for every
 integer $r > 1$ and $i=1,2$
 there are canonical isomorphisms of $\Z_p[G_w]$-modules
 \[
  K_{2r-i}(\mathcal{O}_w; \Z_p) \simeq H^i_{\et}(L_w, \Z_p(r)).
 \]
\end{theorem}

\subsection{Local Galois cohomology}
We keep the notation of \S \ref{sec:setup}. In particular, $L/K$ is a Galois extension
of number fields with Galois group $G$. Let $p$ be an odd prime.
We denote the (finite) set of places of $K$ that ramify in $L/K$ by $S_{\ram}$ and let $S$
be a finite set of places of $K$ containing $S_{\ram}$ and all archimedean and $p$-adic places
(i.e.~$S_{\infty} \cup S_p \cup S_{\ram} \subseteq S$).

Let $M$ be a topological $G_{L,S}$-module.
Then $M$ becomes a topological $G_{L_w}$-module
for every $w \in S(L)$ by restriction. 
For any $i \in \Z$ we put
\[
 P^i(\mathcal O_{L,S}, M) := \bigoplus_{w \in S(L)} H^i_{\et}(L_w, M).
\]

For any integers $r$ and $i$ we define $P^i(\mathcal O_{L,S}, \Q_p(r))$ to be
$P^i(\mathcal O_{L,S}, \Z_p(r)) \otimes_{\Z_p} \Q_p$.

\begin{lemma} \label{lem:local-cohomology-Qp}
 Let $r>1$ be an integer. Then we have isomorphisms of $\Q_p[G]$-modules
 \[
  P^i(\mathcal O_{L,S}, \Q_p(r)) \simeq \left\{
    \begin{array}{ll}
     H_{-r}^+(L) \otimes \Q_p & \mbox{ if } i=0\\
     L \otimes_{\Q} \Q_p & \mbox{ if } i=1\\
     0 & \mbox{ otherwise.}
    \end{array}
  \right.
 \]
\end{lemma}

\begin{proof}
 This is \cite[Lemma 3.3]{wild-kernels} (see also \cite[Lemma 5.2.4]{barrett}).
 The case $i=1$ will be crucial in the following 
 so that we briefly recall its proof.
 Let $w \in S_p(L)$ and put
 $D_{dR}^{L_w}(\Q_p(r)) := H^0(L_w, B_{dR} \otimes_{\Q_p} \Q_p(r))$,
 where $B_{dR}$ denotes Fontaine's de Rham period ring.
 Then the Bloch--Kato exponential map
 \begin{equation} \label{eqn:BK-exponential}
  \exp_r^{BK}: L_w = D_{dR}^{L_w}(\Q_p(r)) \longrightarrow H^1_{\et}(L_w, \Q_p(r))
 \end{equation}
 is an isomorphism for every $w \in S_p(L)$ as follows from \cite[Corollary 3.8.4 and Example 3.9]{MR1086888}. 
 Since the groups $H^1_{\et}(L_w, \Z_p(r))$ are finite for $w \not\in S_p(L)$, we 
 obtain isomorphisms of $\Q_p[G]$-modules
 \[
  P^1(\mathcal O_{L,S}, \Q_p(r)) \simeq \bigoplus_{w \in S_p(L)} H^1_{\et}(L_w, \Q_p(r))
  \simeq \bigoplus_{w \in S_p(L)} L_w \simeq L \otimes_{\Q} \Q_p.
 \]
\end{proof} 

\subsection{Schneider's conjecture}
Let $M$ be a topological $G_{L,S}$-module.
 For any integer $i$
 we denote the kernel of the natural localization map
 \[
  H^i_{\et} (\mathcal O_{L,S}, M) \longrightarrow P^i(\mathcal O_{L,S}, M)
 \]
 by $\Sha^i(\mathcal O_{L,S}, M)$.
 We call $\Sha^i(\mathcal O_{L,S}, M)$
  the \emph{Tate--Shafarevich group} of $M$ in degree $i$.
We recall the following conjecture of Schneider \cite[p.~192]{MR544704}.

\begin{conj}[$\Sch(L,p,r)$] \label{conj:Schneider}
	Let $r \not=0$ be an integer. Then the Tate--Shafarevich group $\Sha^1(\mathcal{O}_{L,S}, \Z_p(r))$ vanishes.
\end{conj}

\begin{remark}
It is not hard to show that Conjecture \ref{conj:Schneider} does not depend on the choice of the set $S$.
\end{remark}

\begin{remark}
 Schneider originally conjectured that $H^2_{\et}(\mathcal O_{L,S}, \Q_p / \Z_p (1-r))$ vanishes. Both conjectures are in fact equivalent 
 (see \cite[Proposition 3.8 (ii)]{wild-kernels}).
\end{remark}

\begin{remark}
	It can be shown that Schneider's conjecture for $r=1$ is equivalent to Leopoldt's conjecture
	(see \cite[Chapter X, \S 3]{MR2392026}).
\end{remark}

\begin{remark} \label{rmk:Schneider-almost}
	For a given number field $L$ and a fixed prime $p$, 
	Schneider's conjecture holds for almost all $r$.
	This follows from \cite[\S 5, Corollar 4]{MR544704} and \cite[\S 6, Satz 3]{MR544704}.
\end{remark}

\begin{remark}
Schneider's conjecture $\Sch(L,p,r)$ holds whenever $r<0$ by work of Soul\'{e}
\cite{MR553999}; see also \cite[Theorem 10.3.27]{MR2392026}.
\end{remark}

\begin{lemma} \label{lem:Sha-torsion-free}
Let $r\not=0$ be an integer and let $p$ be an odd prime.
Then the Tate--Shafarevich group $\Sha^1(\mathcal{O}_{L,S}, \Z_p(r))$ is torsion-free.
In particular, Schneider's conjecture $\Sch(L,p,r)$ holds if and only if
$\Sha^1(\mathcal{O}_{L,S}, \Q_p(r)) := \Q_p \otimes_{\Z_{p}} 
\Sha^1(\mathcal{O}_{L,S}, \Z_p(r))$
vanishes.
\end{lemma}

\begin{proof}
The first claim is \cite[Proposition 3.8 (i)]{wild-kernels}. The second claim is
immediate.
\end{proof}

\subsection{Artin $L$-series}
Let $L/K$ be a finite Galois extension of number fields with Galois group $G$
and let $S$ be a finite set of places of $K$ containing all archimedean places.
For any irreducible complex-valued character $\chi$ of $G$ we denote the $S$-truncated Artin $L$-series
by $L_S(s, \chi)$, and the leading coefficient of $L_S(s, \chi)$ at an integer $r$ by
$L_S^{\ast}(r, \chi)$.
We shall sometimes use this notion even if $L_S^{\ast}(r, \chi) = L_S(r, \chi)$
(which will happen frequently in the following).

Recall that there is a canonical
isomorphism $\zeta(\C[G]) \simeq \prod_{\chi \in \Irr_{\C}(G)} \C$.
We define the equivariant $S$-truncated Artin $L$-series to be the
meromorphic $\zeta(\C[G])$-valued function
\[
 L_S(s) := (L_S(s,\chi))_{\chi \in \Irr_{\C}(G)}.
\]
For any $r \in \Z$ we also put
\[
 L_S^{\ast}(r) := (L_S^{\ast}(r,\chi))_{\chi \in \Irr_{\C}(G)} \in \zeta(\R[G])^{\times}.
\]

\subsection{A conjecture of Gross} \label{subsec:Gross}
Let $r>1$ be an integer.
Since Borel's regulator map \eqref{eqn:regulator} induces an isomorphism of $\R[G]$-modules,
the Noether--Deuring Theorem (see \cite[Lemma 8.7.1]{MR2392026} for instance) implies the existence of $\Q[G]$-isomorphisms
\begin{equation} \label{eqn:phi_1-r}
 \phi_{1-r}: H_{1-r}^+(L) \otimes \Q \stackrel{\simeq}{\longrightarrow} K_{2r-1}(\mathcal{O}_L) \otimes \Q.
\end{equation}
Let $\chi$ be a complex character of $G$ and let $V_{\chi}$ be a $\C[G]$-module with character $\chi$.
Composition with $\rho_r \circ \phi_{1-r}$
induces an automorphism of $\Hom_G(V_{\check\chi}, H_{1-r}^+(L) \otimes \C)$.
Let $R_{\phi_{1-r}}(\chi)  = 
\det_{\C}((\rho_r \circ \phi_{1-r})^{\check \chi}) \in \C^{\times}$ 
be its determinant. If $\chi'$ is a second character,
then $R_{\phi_{1-r}}(\chi + \chi') = 
R_{\phi_{1-r}}(\chi) \cdot R_{\phi_{1-r}}(\chi')$ by Lemma 
\ref{lem:properties-of-rho-parts-functor} so
that we obtain a map
\begin{eqnarray*}
 R_{\phi_{1-r}}: R(G) & \longrightarrow & \C^{\times} \\
 \chi & \mapsto & \det\nolimits_{\C}(\rho_r \circ \phi_{1-r} \mid \Hom_G(V_{\check\chi}, H_{1-r}^+(L) \otimes \C)),
\end{eqnarray*}
where $R(G) := R_{\C}(G)$ denotes the ring of virtual complex characters of $G$.
We likewise define
\begin{eqnarray*}
 A_{\phi_{1-r}}^S: R(G) & \longrightarrow & \C^{\times} \\
 \chi & \mapsto & R_{\phi_{1-r}}(\chi) / L_S^{\ast}(1-r, \chi).
\end{eqnarray*}
Gross \cite[Conjecture 3.11]{MR2154331} conjectured the following higher analogue of Stark's conjecture.

\begin{conj}[Gross] \label{conj:Gross}
 We have $A_{\phi_{1-r}}^S(\chi^{\sigma}) = A_{\phi_{1-r}}^S(\chi)^{\sigma}$ for all $\sigma \in \Aut(\C)$.
\end{conj}

\begin{remark}\label{rmk:Gross-independent-of-choices}
It is not hard to see that Gross's conjecture does not depend on $S$ and the choice of $\phi_{1-r}$
(see also \cite[Remark 6]{MR2801311}). A straightforward substitution shows that if it is true for $\chi$ then it is true for $\chi^{\tau}$ for every choice of $\tau \in \Aut(\C)$.
\end{remark}

We record some cases where Gross's conjecture is known and deduce a few new cases.
If $q = \ell^n$ is a prime power, we let $\Aff(q)$ be the group of affine transformations
on $\F_q$. Thus we may write $\Aff(q) \simeq N \rtimes H$, where 
$H = \{ x \mapsto ax \mid a \in \F_{q}^{\times }\} \simeq \F_q^{\times}$ acts on
$N =\{ x \mapsto x+b \mid b \in \F_{q} \} \simeq \F_q$ in the natural way.
Note that $N$ is the commutator subgroup of $\Aff(q)$.

\begin{theorem} \label{thm:Gross-cases}
Let $L/K$ be a finite Galois extension of number fields with Galois group $G$
and let $\chi \in R(G)$ be a virtual character. Let $r>1$ be an integer.
Then Gross's conjecture (Conjecture \ref{conj:Gross}) holds in each of the following cases.
\begin{enumerate}
\item 
$\chi$ is absolutely abelian, i.e.~there is a normal subgroup $N$ of $G$ such that
$\chi$ factors through $G/N \simeq \Gal(L^N/K)$ and $L^N/\Q$ is abelian;
\item 
$\chi = \mathbbm{1}_{G}$ is the trivial character;
\item 
$\chi$ is a virtual permutation character, i.e.~a $\Z$-linear combination of characters
of the form $\ind_{H}^G \mathbbm{1}_H$ where $H$ ranges over subgroups of $G$;
\item 
$G \simeq \Aff(q) = N \rtimes H$ and $L^N/\Q$ is abelian;
\item
$L^{\ker(\chi)}$ is totally real and $r$ is even;
\item
$L^{\ker(\chi)} / K$ is a CM-extension, $\chi$ is an odd character and $r$ is odd.
\end{enumerate}
\end{theorem}

\begin{proof}
We first note that (ii) is Borel's result \eqref{eqn:Borels-rationality} above.
Since Gross's conjecture is invariant under induction and respects addition of characters,
(ii) implies (iii). For (i), (v) and (vi) we refer the reader to
\cite[Theorem 5.2]{wild-kernels} and the references given therein.
We now prove (iv). It suffices to show that Gross's conjecture holds for every
$\chi \in \Irr_{\C}(G)$. If $\chi$ is linear, it factors through $G/N$
so that $\chi$ is indeed absolutely linear. Thus Gross's conjecture holds 
by (i). It has been shown in the proof of \cite[Theorem 10.5]{MR4111943}
that there is a unique non-linear irreducible character $\chi_{\nl}$ of $G$ 
and that this character
can be expressed as a $\Z$-linear combination of $\ind_H^G \mathbbm{1}_H$
and linear characters in $\Irr_{\C}(G)$. As Gross's conjecture holds for
the linear characters and for $\ind_H^G \mathbbm{1}_H$ by (iii),
it also holds for $\chi_{\nl}$.
\end{proof}

For any integer $k$ we write
\begin{equation} \label{eqn:iota_k}
  \iota_k: L \otimes_{\Q} \C \longrightarrow \bigoplus_{\Sigma(L)} \C =
  \left(H_{1-k}^+(L) \oplus H_{-k}^+(L)\right) \otimes \C
\end{equation}
for the canonical $\C[G \times \Gal(\C / \R)]$-equivariant isomorphism which is induced by mapping
$l \otimes z$ to $(\sigma(l)z)_{\sigma \in \Sigma(L)}$ for $l \in L$
and $z \in \C$. Now fix an integer $r>1$.
We define a $\C[G]$-isomorphism
\begin{eqnarray} \label{eqn:lambda_r}
 \lambda_r: \left(K_{2r-1}(\mathcal{O}_L) \oplus H_{-r}^+(L)\right) \otimes \C & \simeq &
    \left(H_{1-r}^+(L) \oplus H_{-r}^+(L)\right) \otimes \C \\
    & \simeq & L \otimes_{\Q} \C. \nonumber
\end{eqnarray}
Here, the first isomorphism is induced by $\rho_r \oplus \id_{H_{-r}^+(L)}$,
whereas the second isomorphism is $\iota_r^{-1}$. 
As above, there exist $\Q[G]$-isomorphisms
\begin{equation} \label{eqn:phi_r}
 \phi_r: L \stackrel{\simeq}{\longrightarrow} \left(K_{2r-1}(\mathcal{O}_L) \oplus H_{-r}^+(L)\right) \otimes \Q.
\end{equation}
We now define maps
\begin{eqnarray*}
 R_{\phi_{r}}: R(G) & \longrightarrow & \C^{\times} \\
 \chi & \mapsto & \det\nolimits_{\C}\left(\lambda_r \circ \phi_{r} \mid
  \Hom_G(V_{\check\chi}, L \otimes_{\Q} \C)\right) 
\end{eqnarray*}
and
\begin{eqnarray*}
 A_{\phi_{r}}^S: R(G) & \longrightarrow & \C^{\times} \\
 \chi & \mapsto & R_{\phi_{r}}(\chi) / L_S(r, \check\chi).
\end{eqnarray*}

\begin{prop} \label{prop:Gross-equivalent}
 Fix an integer $r>1$ and a character $\chi$. Then Gross's conjecture \ref{conj:Gross}
 holds if and only if we have $A_{\phi_{r}}^S(\chi^{\sigma}) = A_{\phi_{r}}^S(\chi)^{\sigma}$ 
 for all $\sigma \in \Aut(\C)$.
\end{prop}

\begin{proof}
 This is \cite[Proposition 5.5]{wild-kernels}.
\end{proof}

\subsection{The comparison period} \label{sec:comparison}
We henceforth assume that $p$ is an odd prime and that $L/K$ is a Galois CM-extension.
Recall that $\tau \in G$ is the unique automorphism induced by complex conjugation.
For each $n \in \Z$ we define a central idempotent 
$e_n := \frac{1- (-1)^n \tau}{2}$ in $\Z[\half][G]$. 
Now let $r>1$ be an integer.
Since $L$ is CM, the idempotent $e_r$ acts trivially on
$H_{1-r}^+(L) \otimes \C$,
whereas $e_r(H_{-r}^+(L) \otimes \C)$ vanishes. Thus \eqref{eqn:lambda_r} induces a $\C[G]$-isomorphism
\[
 \mu_{\infty}(r): K_{2r-1}(\mathcal{O}_L) \otimes \C \stackrel{\sim}{\longrightarrow} e_r (L \otimes_{\Q} \C).
\]
We likewise define a $\C_p[G]$-homomorphism
\[
    \mu_p(r): K_{2r-1}(\mathcal{O}_L) \otimes \C_p \longrightarrow e_r (L \otimes_{\Q} \C_p)
\]
as follows. For each $w \in S_p(L)$ we 
let $H^i_{\syn}(\mathcal{O}_w, r)$ be the $i$-th syntomic cohomology group as considered
by Besser \cite{MR1809626}. We let $\reg_r^w: K_{2r-1}(\mathcal{O}_w) \rightarrow H^1_{\syn}(\mathcal{O}_w, r)$
be the syntomic regulator \cite[Theorem 7.5]{MR1809626}. By \cite[Lemma 2.15]{MR2520465} (which heavily relies on 
\cite[Proposition 8.6]{MR1809626}) we have 
canonical isomorphisms $H^1_{\syn}(\mathcal{O}_w, r) \simeq L_w$
for each $w \in S_p(L)$. The map $\mu_p(r)$ is induced by the following chain of homomorphisms
\begin{eqnarray} \label{eqn:maps-for-mu_p(r)}
 K_{2r-1}(\mathcal{O}_L) & \longrightarrow & \bigoplus_{w \in S_p(L)} K_{2r-1}(\mathcal{O}_w)\\
 & \longrightarrow & \bigoplus_{w \in S_p(L)} H^1_{\syn}(\mathcal{O}_w, r)\nonumber \\
 & \simeq & \bigoplus_{w \in S_p(L)} L_w\nonumber \\
 & \simeq & L \otimes_{\Q} \Q_p. \nonumber 
\end{eqnarray}
The map $\mu_p(r)$ shows up in the formulation of the $p$-adic Beilinson conjecture.
However, the following map will be more suitable for the relation to the ETNC.
We define a $\C_p[G]$-homomorphism
\begin{eqnarray*}
 \tilde \mu_p(r): K_{2r-1}(\mathcal{O}_L) \otimes \C_p & = & K_{2r-1}(\mathcal{O}_{L,S}) \otimes \C_p\\
 & \simeq & H^1_{\et}(\mathcal{O}_{L,S}, \Z_p(r)) \otimes_{\Z_p} \C_p\\
 & \longrightarrow & e_r(P^1(\mathcal{O}_{L,S}, \Z_p(r)) \otimes_{\Z_p} \C_p)\\
 & \simeq & e_r (L \otimes_{\Q} \C_p).
\end{eqnarray*}
Here, the first map is induced by the $p$-adic Chern class map $ch_{r,1}^{(p)}$
which is an isomorphism by Theorem \ref{thm:Quillen--Lichtenbaum}; the arrow is the natural
localization map, and the last isomorphism is induced by the Bloch--Kato exponential maps (see
Lemma \ref{lem:local-cohomology-Qp}).

The following result will be crucial for relating the $p$-adic Beilinson
conjecture to the ETNC.

\begin{prop} \label{prop:crucial}
 For each $r>1$ we have $\mu_p(r) = \tilde\mu_p(r)$.
\end{prop}

\begin{proof}
 For any abelian group
 $A$ we write $\widehat A$ for its $p$-completion, that is
 $\widehat A := \varprojlim_n A / p^n A$. 
 The localization maps \eqref{eqn:maps-for-mu_p(r)} induce
 a map 
 \[
  K_{2r-1}(\mathcal{O}_L) \otimes \Z_p \longrightarrow \bigoplus_{w \in S_p(L)} \widehat{K_{2r-1}(\mathcal{O}_w)}.
 \]
 For each $w \in S_p(L)$ the Universal Coefficient Theorem \cite[Chapter IV, Theorem 2.5]{MR3076731}
 implies that there is a natural (injective) map 
 \[
   \widehat{K_{2r-1}(\mathcal{O}_w)} \longrightarrow K_{2r-1}(\mathcal{O}_w; \Z_p).
 \]
 Moreover, by \cite[Corollary 9.10]{MR1809626} there is a natural map 
 $H^1_{\syn}(\mathcal{O}_w, r) \rightarrow H^1_{\et}(L_w, \Q_p(r))$ such that the diagram
 \[ \xymatrix{
   \widehat{K_{2r-1}(\mathcal{O}_w)} \ar[r] \ar[d]  & K_{2r-1}(\mathcal{O}_w; \Z_p) \ar[d]\\
   H^1_{\syn}(\mathcal{O}_w, r) \ar[r] &  H^1_{\et}(L_w, \Q_p(r))
 }\]
 commutes. Here, the left-hand vertical arrow is induced by the syntomic regulator, and the map on the right by the isomorphism in 
 Theorem \ref{thm:local-Quillen--Lichtenbaum}. Moreover,
 the composite map 
 \[
  L_w \simeq H^1_{\syn}(\mathcal{O}_w, r) \longrightarrow H^1_{\et}(L_w, \Q_p(r))
 \]
 is the Bloch--Kato exponential map \eqref{eqn:BK-exponential} 
 by \cite[Proposition 9.11]{MR1809626}.
 Unravelling the definitions we now see that the maps $\mu_p(r)$ and $\tilde\mu_p(r)$
 coincide.
\end{proof}

We will henceforth often not distinguish between the maps $\mu_p(r)$ and $\tilde\mu_p(r)$.
Since the Tate--Shafarevich group $\Sha^1(\mathcal{O}_{L,S}, \Z_p(r))$ is torsion-free by Lemma \ref{lem:Sha-torsion-free},
the following result is now immediate.

\begin{lemma} \label{lem:mu_p-Schneider}
 The map $\mu_p(r)$ is a $\C_p[G]$-isomorphism if and only if $\Sch(L,p,r)$ holds.
\end{lemma}

\begin{definition}
Let $j: \C \simeq \C_{p}$ be a field isomorphism and let $\rho \in R_{\C_{p}}^{+}(G)$. Let $r>1$ be an integer.
We define the comparison period attached to $j$, $\rho$ and $r$ to be
\[ 
\Omega_{j}(r, \rho) := 
{\det}_{\C_{p}}(\mu_{p}(r) \circ (\C_{p} \otimes_{\C,j} \mu_{\infty}(r))^{-1})^{\rho} \in \C_{p}.
\] 
\end{definition}

We record some basic properties of $\Omega_{j}(r,-)$.

\begin{lemma}\label{lem:properties-of-omega}
Let $H,N$ be subgroups of $G$ with $N$ normal in $G$.
\begin{enumerate}
\item 
Let $\rho_{1}, \rho_{2} \in  R_{\C_{p}}^{+}(G)$. Then $\Omega_{j}(r, \rho_{1} + \rho_{2}) = \Omega_{j}(r, \rho_{1}) \Omega_{j}(r, \rho_{2})$.
\item
Let $\rho \in R_{\C_{p}}^{+}(H)$. Then $\Omega_{j}(r, \ind_{H}^{G} \rho) = \Omega_{j}(r,\rho)$.
\item
Let $\rho \in R_{\C_{p}}^{+}(G / N)$. Then $\Omega_{j}(r, \infl_{G / N}^{G} \rho) = \Omega_{j}(r,\rho)$.
\end{enumerate}
\end{lemma}

\begin{proof}
Each part follows from the corresponding part of Lemma \ref{lem:properties-of-rho-parts-functor}.
\end{proof}

\begin{remark}\label{rmk:vanishing-of-mu-p-independent-of-j}
Since $\mu_{\infty}(r)$ is an isomorphism, for any two choices of field isomorphism $j,j': \C \simeq \C_{p}$ we have that 
$\Omega_{j}(r,\rho)=0$ if and only if $\Omega_{j'}(r,\rho)=0$.
\end{remark}

\begin{remark}\label{rmk:rho-parts-of-Leo}
For any fixed choice of field isomorphism $j:\C \simeq \C_{p}$ we have 
\begin{eqnarray*}
\Sch(L,p,r) \textrm{ holds} 
&\Longleftrightarrow& \mu_{p}(r) \textrm{ is an isomorphism} \\
&\Longleftrightarrow& \Omega_{j}(r,\rho) \neq 0 \quad \forall \rho \in \Irr_{\C_{p}}(G),
\end{eqnarray*}
where the first equivalence is Lemma \ref{lem:mu_p-Schneider}.
Thus the non-vanishing of $\Omega_{j}(r,\rho)$ can be thought of as the `$\rho$-part' of $\Sch(L,p,r)$.
Moreover, if $\Omega_{j}(r,\rho) \neq 0$ then we may set $\Omega_{j}(r, -\rho) := \Omega_{j}(r,\rho)^{-1}$
and so if we assume $\Sch(L,p,r)$ then Lemma \ref{lem:properties-of-omega} (i) shows that the definition of $\Omega_{j}(r,\rho)$ 
naturally extends to any virtual character $\rho \in R_{\C_{p}}(G)$.
\end{remark}

\begin{remark}
	Assume that $G$ is abelian. For each integer $r$,
	Burns, Kurihara and Sano \cite[\S 2.2]{BKS-Stark} define canonical
	period-regulator isomorphisms
	\begin{equation} \label{eqn:BKS-regulator}
		\varepsilon_r \bigwedge_{\C_p[G]}^{r_j^{\varepsilon}} H^1_{\et}(\mathcal{O}_{L,S}, \Z_p(1-r)) \otimes_{\Z_p} \C_p \rightarrow \varepsilon_r 
		\bigwedge_{\C_p[G]}^{r_j^{\varepsilon}}
		P^0(\mathcal{O}_{L,S}, \Z_p(-r)) \otimes_{\Z_p} \C_p.
	\end{equation}
	Here $\varepsilon_r \in \Z_p[G]$ are certain idempotents such that
	the $\varepsilon_r$-parts of both $\C_p[G]$-modules 
	in the exterior products
	are free of the same rank $r_j^{\varepsilon}$.
	If $r>1$ and $\Sch(L,p,r)$ holds, then one may take 
	$\varepsilon_r = e_r$ and $r_j^{\epsilon} = 0$.
	In this case the diagram \cite[p.~125]{BKS-Stark} gives an exact
	sequence of $\C_p[G]$-modules
	\[
	0 = e_r (H^1_{\et}(\mathcal{O}_{L,S}, \Z_p(1-r)) \otimes_{\Z_p} \C_p)
	\rightarrow	e_r(L \otimes_{\Q} \C_p)^{\ast} \rightarrow 
	e_r(H_{1-r}^+(L) \otimes \C_p)^{\ast} \rightarrow 0,
	\]
	where $(-)^{\ast}$ denotes $\C_p$-linear duals (note also that
	our $H_{1-r}^+(L)$ is their $H_L(r-1)^+$). The non-trivial
	map is (up to sign) the dual of $\mu_p(r) \circ \rho_r^{-1}$.
	Hence the exterior product on the left of \eqref{eqn:BKS-regulator}
	canonically identifies with 
	\[
		e_r \det\nolimits_{\C_p[G]}((L \otimes_{\Q} \C_p)^{\ast})
		\otimes_{\C_p[G]} \det\nolimits_{\C_p[G]}^{-1}((H_{1-r}^+(L) \otimes \C_p)^{\ast})
	\]
	and the isomorphism $\mu_p(r) \circ \rho_r^{-1}$ together with $\iota_r$
	induces a map to $\det_{\C_{p}[G]}(H_{-r}^+(L) \otimes \C_p)$ which
	can be identified with the exterior power on the right by a variant of Lemma
	\ref{lem:local-cohomology-Qp}.
	For more details we refer the
	interested reader to \cite[\S 2.2.4]{BKS-Stark}.
	
	The authors then use the isomorphism \eqref{eqn:BKS-regulator}  
	to define generalized Stark elements and to state 
	\cite[Conjecture 3.6]{BKS-Stark} which might be seen as an analogue
	and refinement of a conjecture of Rubin \cite{MR1385509} in the case $r=0$.
	It is then shown in \cite[\S 4]{BKS-Stark} that their conjecture
	is implied by the appropriate special case of the ETNC. As the formulation
	of the latter involves the Bloch--Kato exponential map rather than the
	syntomic regulator, a variant of Proposition \ref{prop:crucial}
	is already implicit in their work 
	(for instance, see \cite[Remark 2.7]{BKS-Stark})
\end{remark}

\subsection{$p$-adic Artin $L$-functions}\label{subsec:p-adic-Artin-L-functions}
Let $E/K$ be a finite Galois extension of totally real number fields and let $G=\Gal(E/K)$. Let $p$ be an odd prime
and let $S$ be a finite set of places of $K$ containing $S_{p} \cup S_{\infty}$.
For each $\rho \in R_{\C_{p}}(G)$ the $S$-truncated $p$-adic Artin $L$-function 
attached to $\rho$ is the unique $p$-adic meromorphic
function $L_{p,S}(s, \rho) : \Z_{p} \rightarrow \C_{p}$ with the property that for each strictly negative 
integer $n$ and each field isomorphism $j: \C \simeq \C_{p}$ we have
\[
L_{p,S}(n, \rho) = j \left( L_{S}(n, (\rho \otimes \omega^{n-1})^{j^{-1}}) \right),
\]
where $\omega: G_K \rightarrow \Z_{p}^{\times}$ is the Teichm\"uller character
and we view $\rho \otimes \omega^{n-1}$ as a character of $\Gal(E(\zeta_p)/K)$.
By a result of Siegel \cite{MR0285488} the right-hand side 
does indeed not depend on the choice of $j$.
In the case that $\rho$ is linear, $L_{p,S}(s, \rho)$ was constructed independently by 
Deligne and Ribet \cite{MR579702}, Barsky \cite{MR525346} and Cassou-Nogu\'es \cite{MR524276}.
Greenberg \cite{MR692344} then extended the construction to the general case using Brauer induction.

\subsection{Statement of the $p$-adic Beilinson conjecture}
We now formulate our variant of the $p$-adic Beilinson conjecture.

\begin{conj}[The $p$-adic Beilinson conjecture]\label{conj:p-adic-Beilinson}
Let $E/K$ be a finite Galois extension of totally real number fields and let $G=\Gal(E/K)$.
Let $p$ be an odd prime and let $S$ be a finite set of places of $K$ containing 
$S_{p} \cup S_{\infty}$.
Let $\rho \in R_{\C_{p}}^{+}(G)$ and let $r>1$ be an integer. 
Then for every choice of field isomorphism $j: \C \cong \C_{p}$ we have 
\begin{equation}\label{eqn:p-adic-Beilinson-at-r}
L_{p, S}(r,\rho) = \Omega_{j}(r, \rho \otimes \omega^{r-1}) \cdot 
 j \left( L_{S}(r, (\rho \otimes \omega^{r-1})^{j^{-1}}) \right). 
\end{equation}
\end{conj}

\begin{remark}
 It is straightforward to show that Conjecture \ref{conj:p-adic-Beilinson} does not depend on
 the choice of $S$.
\end{remark}

\begin{remark} \label{rem:Schneider-vs-L_p}
 One can show (see Theorem \ref{thm:Schneider-semisimple} below) 
 that $L_{p, S}(r,\rho) \not= 0$ if and only if $\Omega_{j}(r, \rho \otimes \omega^{r-1})
 \not= 0$. In this case (and thus in particular if $\Sch(E(\zeta_p),p,r)$ holds) 
 the statement of Conjecture \ref{conj:p-adic-Beilinson} naturally extends to all virtual
 characters $\rho \in R_{\C_p}(G)$.
\end{remark}

\begin{remark} \label{rem:compatible-Besser-etal}
 It is clear from the definitions that Conjecture \ref{conj:p-adic-Beilinson}
 is compatible with the $p$-adic Beilinson conjecture as considered by
 Besser, Buckingham, de Jeu and Roblot \cite[Conjecture 3.18]{MR2520465}.
 More concretely, the equality \eqref{eqn:p-adic-Beilinson-at-r} 
 is equivalent to the appropriate special case of \cite[Conjecture 3.18 (i)--(iii)]{MR2520465},
 whereas \cite[Conjecture 3.18 (iv)]{MR2520465} then is equivalent to the non-vanishing of
 $\Omega_{j}(r, \rho \otimes \omega^{r-1})$ as follows from Remark \ref{rem:Schneider-vs-L_p}.
\end{remark}

\begin{remark} \label{rem:infl-ind-invariance}
Since both complex and $p$-adic Artin $L$-functions satisfy properties analogous to those 
of $\Omega_{j}(r, -)$ given in Lemma \ref{lem:properties-of-omega}, 
the truth of Conjecture \ref{conj:p-adic-Beilinson} is invariant under induction and inflation;
moreover, if it holds for $\rho_{1}, \rho_{2} \in R_{\C_{p}}(G)$ then it holds for $\rho_{1}+\rho_{2}$.
\end{remark}

\subsection{The relation to Gross's conjecture} \label{subsec:Gross-relation}

The following results are the analogues of \cite[Theorem 4.16 and Corollary 4.18]{MR4111943}, respectively.

\begin{theorem}\label{thm:p-Beilinson-implies-Gross}
Let $E/K$ be a finite Galois extension of totally real number fields and let $G=\Gal(E/K)$.
Let $p$ be an odd prime and let $S$ be a finite set of places of $K$ containing 
$S_{p} \cup S_{\infty}$.
Let $\rho \in R_{\C_{p}}^{+}(G)$ and let $r>1$ be an integer. We put $\psi := \rho \otimes \omega^{r-1}$.
If $\Omega_{j}(r,\psi) \neq 0$ for some (and hence every) choice of field isomorphism $j:\C \cong \C_{p}$ 
then the following statements are equivalent.
\begin{enumerate}
\item $\Omega_{j}(r,\psi) \cdot j ( L_{S}(r, \psi^{j^{-1}}) )$ is independent of the choice of $j:\C \cong \C_{p}$.
\item Gross's conjecture at $s=1-r$ holds for $\check\psi^{j^{-1}} \in R_{\C}^{+}(\Gal(E(\zeta_p)/K)$ and some (and hence every) choice of $j:\C \cong \C_{p}$.
\end{enumerate}
\end{theorem}

\begin{proof}
The first and second occurrence of `and hence every' in the 
statement of the theorem follow from Remark 
\ref{rmk:Gross-independent-of-choices} and Remark \ref{rmk:vanishing-of-mu-p-independent-of-j}, respectively.

Let $j,j': \C \cong \C_{p}$ be field isomorphisms and let $\chi := \psi^{j^{-1}}$.
Then $j = j' \circ \sigma$ for some $\sigma \in \Aut(\C)$ and so $\psi^{j'^{-1}} = \chi^{\sigma}$.
For every $\Q[G]$-isomorphism
$\phi_r$ as in \eqref{eqn:phi_r} we have
\[
\Omega_{j}(r,\psi) \cdot j \left( R_{\phi_r}(\check\chi) \right) = {\det}_{\C_{p}} (\mu_{p}(r) \circ (\C_{p} \otimes_{\Q} \phi_r))^{\psi},
\]
which does not depend on $j$. 
Hence we have
\[
\frac{\Omega_{j}(r,\psi) \cdot j ( L_{S}(r, \psi^{j^{-1}}) )}{\Omega_{j'}(r,\psi) \cdot  j' ( L_{S}(r, \psi^{(j')^{-1}}) )}
= \frac{j'(R_{\phi_r}(\check\chi)) \cdot j(L_{S}(r, \chi))}{j(R_{\phi_r}(\check\chi)) \cdot j'(L_{S}(r, \chi^{\sigma}))}
= j' \left( \frac{\sigma(L_{S}(r, \chi)) \cdot R_{\phi_r}(\check\chi)}
{\sigma(R_{\phi_r}(\check\chi)) \cdot L_{S}(r, \chi^{\sigma})} \right).
\]
By Proposition \ref{prop:Gross-equivalent} the last expression is equal to $1$ if and only if Gross's conjecture at $s=1-r$ (Conjecture \ref{conj:Gross}) holds for the character $\check\chi$.
\end{proof}

\begin{corollary}\label{cor:p-Beilinson-implies-Gross}
Let $E/K$ be a finite Galois extension of totally real number fields and let $G=\Gal(E/K)$.
Fix a prime $p$ and let $r>1$ be an integer. Assume that $\Sch(E(\zeta_p),p,r)$ holds.
If the $p$-adic Beilinson conjecture at $s=r$ holds for all $\rho \in R_{\C_{p}}^{+}(G)$ 
then Gross's conjecture at $s=1-r$ holds for $\chi \otimes \omega^{r-1}$ for all $\chi \in R_{\C}^{+}(G)$.  
\end{corollary} 

\subsection{Absolutely abelian characters} \label{subsec:absolutely}
Since our conjecture is compatible with that of \cite{MR2520465} by 
Remark \ref{rem:compatible-Besser-etal} and invariant under induction and inflation of characters
by Remark \ref{rem:infl-ind-invariance}, we deduce the following result from work of Coleman \cite{MR674400}
(see \cite[Proposition 4.17]{MR2520465}).

\begin{theorem} \label{thm:p-Beilinson-abelian}
 Let $E/K$ be a finite Galois extension of totally real number fields and let $G=\Gal(E/K)$. 
 Let $p$ be a prime and let $r>1$ be an integer.
Suppose that $\rho \in R_{\C_{p}}^+(G)$ is an absolutely abelian character, i.e.,
there exists a normal subgroup $N$ of $G$ such that $\rho$ factors through $G/N \cong \Gal(E^{N}/K)$ and $E^{N}/\Q$ is abelian. 
Then the $p$-adic Beilinson conjecture (Conjecture \ref{conj:p-adic-Beilinson}) holds for $\rho$.
\end{theorem}

\section{Equivariant Iwasawa theory}

\subsection{Bockstein homomorphisms}
We recall some background material regarding Bockstein homomorphisms. The reader may also
consult \cite[\S 3.1-\S 3.2]{MR2290587}. 

Let $\mathcal{G}$ be a compact $p$-adic Lie group that contains a closed normal 
subgroup $H$ such that $\Gamma := \mathcal{G}/H$ is isomorphic to $\Z_p$.
We fix a topological generator $\gamma$ of $\Gamma$.
The Iwasawa algebra of $\mathcal{G}$ is
\[
\Lambda(\mathcal{G}) := \Z_{p}\llbracket\mathcal{G}\rrbracket = \varprojlim \Z_{p}[\mathcal{G}/\mathcal{N}],
\]
where the inverse limit is taken over all open normal subgroups $\mathcal{N}$ of $\mathcal{G}$.
If $F$ is a finite field extension of $\Q_{p}$  with ring of integers $\mathcal{O}=\mathcal{O}_{F}$,
we put $\Lambda^{\mathcal{O}}(\mathcal{G}) := \mathcal{O} \otimes_{\Z_{p}} \Lambda(\mathcal{G}) = \mathcal{O}\llbracket\mathcal{G}\rrbracket$.
We consider continuous homomorphisms
\begin{equation} \label{eqn:rho}
	\pi: \mathcal{G} \longrightarrow \Aut_{\mathcal{O}}(T_{\pi}),
\end{equation}
where $T_{\pi}$ is a finitely generated free $\mathcal{O}$-module.
For $g \in  \mathcal{G}$ we denote its image in $\Gamma$ under the canonical
projection by $\overline{g}$. We view $\Lambda^{\mathcal{O}}(\Gamma)
\otimes_{\mathcal{O}} T_{\pi}$ as a $(\Lambda^{\mathcal{O}}(\Gamma),
\Lambda(\mathcal{G}))$-bimodule, where $\Lambda^{\mathcal{O}}(\Gamma)$
acts by left multiplication and $\Lambda(\mathcal{G})$ acts on the right via
\[
	(\lambda \otimes_{\mathcal{O}} t) g := \lambda \overline{g} \otimes_{\mathcal{O}} g^{-1}t
\]
for $\lambda \in \Lambda^{\mathcal{O}}(\Gamma)$, $t \in T_{\pi}$ and 
$g \in \mathcal{G}$.
For each complex $C^{\bullet} \in \mathcal{D}^{\perf}(\Lambda(\mathcal{G}))$
we define a complex $C^{\bullet}_{\pi} \in
\mathcal{D}^{\perf}(\Lambda^{\mathcal{O}}(\Gamma))$ by
\[
	C_{\pi}^{\bullet} := (\Lambda^{\mathcal{O}}(\Gamma)
	\otimes_{\mathcal{O}} T_{\pi}) 
	\otimes_{\Lambda(\mathcal{G})}^{\mathbb L} C^{\bullet}.
\]
Given an open normal subgroup $U$ of $\mathcal{G}$ we set
$C_U^{\bullet} := \Z_p[\mathcal{G}/U] 
\otimes_{\Lambda(\mathcal{G})}^{\mathbb L} C^{\bullet}$ and,
if $U$ is contained in the kernel of $\pi$, we furthermore obtain a complex
\[ 
	C_{U,(\pi)}^{\bullet} = 
	T_{\pi} \otimes_{\Z_p[\mathcal{G}/U]}^{\mathbb L} C_U^{\bullet}
	= T_{\pi} \otimes_{\Lambda(\mathcal{G})}^{\mathbb L} C^{\bullet}
	\in \mathcal{D}^{\perf}(\mathcal{O})
\]
which does actually not depend on $U$. 
The natural exact triangles
\[
	C_{\pi}^{\bullet} \xrightarrow{\gamma - 1} C_{\pi}^{\bullet} 
	\longrightarrow C_{U, (\pi)}^{\bullet}
\]
in $\mathcal{D}(\mathcal{O})$ induce short exact sequences 
of $\mathcal{O}$-modules
\begin{equation} \label{eqn:tautological-ses}
0 \longrightarrow H^i(C_{\pi}^{\bullet})_{\Gamma} 
\stackrel{\alpha^i}{\longrightarrow} H^i(C^{\bullet}_{U,(\pi)})
\stackrel{\tilde \alpha^i}{\longrightarrow} H^{i+1}(C^{\bullet}_{\pi})^{\Gamma} \longrightarrow 0
\end{equation}
for each $i \in \Z$. The \emph{Bockstein homomorphism} in degree $i$ is
defined to be the composite homomorphism
\[
\beta^i_{C^{\bullet},\pi}: H^i(C^{\bullet}_{U,(\pi)})
\stackrel{\tilde \alpha^i}{\longrightarrow} H^{i+1}(C^{\bullet}_{\pi})^{\Gamma} 
\longrightarrow H^{i+1}(C^{\bullet}_{\pi})_{\Gamma} 
\stackrel{\alpha^{i+1}}{\longrightarrow} H^{i+1}(C^{\bullet}_{U,(\pi)}),
\]
where the middle arrow is the tautological map. 
We obtain a bounded complex of $\mathcal{O}$-modules
\[
	\Delta_{C^{\bullet},\pi}: \dots \xrightarrow{\beta^{i-1}_{C^{\bullet},\pi}} 
	H^i(C^{\bullet}_{U,(\pi)})
	\xrightarrow{\beta^{i}_{C^{\bullet},\pi}} H^{i+1}(C^{\bullet}_{U,(\pi)})
	\xrightarrow{\beta^{i+1}_{C^{\bullet},\pi}} \dots
\]
where the term $H^i(C^{\bullet}_{U,(\pi)})$ is placed in degree $i$.
Note that the Bockstein homomorphisms and the complex $\Delta_{C^{\bullet},\pi}$
actually depend on the choice of $\gamma$, though our notation 
does not reflect this. 
The complex $C^{\bullet}$ is called \emph{semisimple at $\pi$} if
$\Q_p \otimes_{\Z_p}\Delta_{C^{\bullet},\pi}$ is acyclic
(for any and hence every choice of $\gamma$).

For any $\Lambda(\mathcal{G})$-module $M$ we put
$M_{\pi} := (\Lambda^{\mathcal{O}}(\Gamma) \otimes_{\mathcal{O}} T_{\pi}) 
\otimes_{\Lambda(\mathcal{G})} M$.

\begin{lemma} \label{lem:one-degree-pi}
Assume that $\mathcal{G}$ is a compact $p$-adic Lie group of dimension $1$ and
let $C^{\bullet} \in \mathcal{D}^{\perf}(\Lambda(\mathcal{G}))$
be acyclic outside degree $a$ for some $a \in \Z$. Further assume
that the $\Lambda(\mathcal{G})$-module $H^a(C^{\bullet})$ 
has projective dimension at most $1$. Then for every continuous
homomorphism $\pi: \mathcal{G} \rightarrow \Aut_{\mathcal{O}}(T_{\pi})$
the complex $C^{\bullet}_{\pi}$ is acyclic outside degree $a$
and there is a canonical isomorphism of $\Lambda^{\mathcal{O}}(\Gamma)$-modules
$H^a(C^{\bullet}_{\pi}) \simeq H^a(C^{\bullet})_{\pi}$.
\end{lemma}

\begin{proof}
This has been shown in the course of the proof of \cite[Lemma 5.6]{MR4092926}.
We repeat the short argument for convenience. 

We may assume that $a=0$. Then $C^{\bullet}$ may be represented by a complex
$P^{-1} \stackrel{d}{\rightarrow} P^0$, where $P^{-1}$ and $P^0$ are projective
$\Lambda(\mathcal{G})$-modules placed in degrees $-1$ and $0$, respectively,
and the homomorphism $d$ is injective.
Then $C^{\bullet}_{\pi}$ is represented by
\[
	(T_{\pi} \otimes_{\Z_{p}} P^{-1})^H \stackrel{d'}{\longrightarrow}
	(T_{\pi} \otimes_{\Z_{p}} P^{0})^H,
\]
where $d'$ is injective since $d$ is. The result follows.
\end{proof}

\subsection{Algebraic $K$-theory}\label{subsec:K-theory}
Let $R$ be a noetherian integral domain with field of fractions $E$.
Let $A$ be a finite-dimensional semisimple $E$-algebra and let $\mathfrak{A}$ be an $R$-order in $A$. For any field extension $F$ of $E$ we set
$A_F := F \otimes_E A$.
Let $K_0(\mathfrak A, F) = K_{0}(\mathfrak{A}, A_F)$ 
denote the relative algebraic $K$-group associated to the ring homomorphism
$\mathfrak{A} \hookrightarrow A_F$.
We recall that $K_{0}(\mathfrak{A}, A_F)$ is an abelian group with generators $[X,g,Y]$ where
$X$ and $Y$ are finitely generated projective $\mathfrak{A}$-modules
and $g:F \otimes_{R} X \rightarrow F \otimes_{R} Y$ is an isomorphism of $A_F$-modules;
for a full description in terms of generators and relations, we refer the reader to \cite[p.\ 215]{MR0245634}.
Moreover, there is a long exact sequence of relative $K$-theory (see \cite[Chapter 15]{MR0245634})
\begin{equation}\label{eqn:long-exact-seq}
K_{1}(\mathfrak{A}) \longrightarrow K_{1}(A_F)
\stackrel{\partial}{\longrightarrow} K_{0}(\mathfrak{A}, A_F)
\longrightarrow K_{0}(\mathfrak{A}) \longrightarrow K_{0}(A_F).
\end{equation}
The reduced norm map $\Nrd = \Nrd_{A_F}: A_F \rightarrow \zeta(A_F)$ is defined componentwise on the Wedderburn decomposition of $A_F$
and extends to matrix rings over $A_F$ (see \cite[\S 7D]{MR632548}); thus
it induces a map $K_{1}(A_F) \longrightarrow \zeta(A_F)^{\times}$, which we also denote by $\Nrd$.

In the case $E = F$ the relative $K$-group $K_{0}(\mathfrak{A}, A)$ identifies with the Grothendieck group whose generators are $[C^{\bullet}]$,
where $C^{\bullet}$
is an object of the category $\mathcal C^{b}_{\tors}(\PMod(\mathfrak{A}))$ of bounded complexes of finitely generated projective $\mathfrak{A}$-modules whose 
cohomology modules are $R$-torsion, and the relations are as follows: $[C^{\bullet}] = 0$ if $C^{\bullet}$ is acyclic, and
$[C_{2}^{\bullet}] = [C_{1}^{\bullet}] + [C_{3}^{\bullet}]$ for every short exact sequence
\[
0 \longrightarrow C_{1}^{\bullet} \longrightarrow C_{2}^{\bullet} \longrightarrow C_{3}^{\bullet} \longrightarrow 0
\]
in $\mathcal C^{b}_{\tors}(\PMod(\mathfrak{A}))$ (see \cite[Chapter 2]{MR3076731} or \cite[\S 2]{MR3068893}, for example).

We denote the full triangulated subcategory of $\mathcal{D}^{\perf} (\mathfrak{A})$
comprising perfect complexes whose cohomology modules are $R$-torsion by 
$\mathcal{D}^{\perf}_{\tors} (\mathfrak{A})$.
Then every object $C^{\bullet}$ of $\mathcal{D}^{\perf}_{\tors} (\mathfrak{A})$ defines a 
class $[C^{\bullet}]$ in $K_{0}(\mathfrak{A}, A)$.

Let $p$ be an odd prime and let $\mathcal{G}$ 
be a one-dimensional compact $p$-adic Lie group that surjects onto $\Z_p$.
Then $\mathcal{G}$ may be written as $\mathcal{G} = H \rtimes \Gamma$
with a finite normal subgroup $H$ of $\mathcal{G}$ and a subgroup
$\Gamma \simeq \Z_p$. Let $\mathcal{O}$ be the ring of integers in some
finite extension $F$ of $\Q_p$.
We consider $\mathfrak{A} = \Lambda^{\mathcal{O}}(\mathcal{G})
= \mathcal O\llbracket\mathcal{G}\rrbracket$ as an order over 
$R := \Lambda^{\mathcal{O}}(\Gamma_0)$,
where $\Gamma_{0}$ is an open subgroup of $\Gamma$ that is central in $\mathcal{G}$. We denote the fraction field of
$\Lambda^{\mathcal{O}}(\Gamma_0)$ by $\mathcal{Q}^F(\Gamma_0)$ and
let $A = \mathcal{Q}^F(\mathcal{G}) :=
\mathcal{Q}^F(\Gamma_0) \otimes_R \Lambda^{\mathcal{O}}(\mathcal{G})$ be the total ring
of fractions of $\Lambda^{\mathcal{O}}(\mathcal{G})$.
Then \cite[Corollary 3.8]{MR3034286}  
shows that the map $\partial$ in \eqref{eqn:long-exact-seq}
is surjective;
thus the sequence
\begin{equation}\label{eqn:Iwasawa-K-sequence}
K_{1}(\Lambda^{\mathcal{O}}(\mathcal{G})) \longrightarrow K_{1}(\mathcal{Q}^F(\mathcal{G})) \stackrel{\partial}{\longrightarrow}
K_{0}(\Lambda^{\mathcal{O}}(\mathcal{G}),\mathcal{Q}^F(\mathcal{G}))
 \longrightarrow 0
\end{equation}
is exact. 
If $\xi \in K_{1}(\mathcal{Q}^F(\mathcal{G}))$ is a pre-image of
some $x \in K_{0}(\Lambda^{\mathcal{O}}(\mathcal{G}),\mathcal{Q}^F(\mathcal{G}))$,
we say that $\xi$ is a \emph{characteristic element} for $x$.
We also set $\mathcal{Q}(\mathcal{G}) := 
\mathcal{Q}^{\Q_p}(\mathcal{G})$.\\

We include the following consequence of \eqref{eqn:Iwasawa-K-sequence} for later use.

\begin{lemma} \label{lem:free-resolution}
	Let $M$ be a finitely generated $\Lambda^{\mathcal{O}}(\mathcal{G})$-module
	of projective dimension at most one. Assume that $M$ is torsion 
	as an $R$-module. Then $M$ admits a free resolution of the form
	\begin{eqnarray} \label{eqn:free-resolution}
		0 \rightarrow \Lambda^{\mathcal{O}}(\mathcal{G})^m \rightarrow
		\Lambda^{\mathcal{O}}(\mathcal{G})^m \rightarrow M \rightarrow 0
	\end{eqnarray}
	for some positive integer $m$. 
	
	Let $\Gamma' \simeq \Z_p$ be an open normal subgroup of $\mathcal{G}$
	and set $G := \mathcal{G}/\Gamma'$.
	If in addition $M_{\Gamma'}$ is finite,
	then $M^{\Gamma'}$ vanishes and \eqref{eqn:free-resolution} induces
	a short exact sequence of $\Z_p[G]$-modules
	\[
	0 \rightarrow \Z_p[G]^m \rightarrow
	\Z_p[G]^m \rightarrow M_{\Gamma'} \rightarrow 0.
	\]
\end{lemma}

\begin{proof}
	Choose a surjective $\Lambda^{\mathcal{O}}(\mathcal{G})$-homomorphism
	$\Lambda^{\mathcal{O}}(\mathcal{G})^m \rightarrow M$. Its kernel is
	a projective $\Lambda^{\mathcal{O}}(\mathcal{G})$-module $P$ by assumption.
	Since $M$ is $R$-torsion, we see that the classes of $P$ and
	$\Lambda^{\mathcal{O}}(\mathcal{G})^m$ in $K_0(\Lambda^{\mathcal{O}}(\mathcal{G}))$
	have the same image in $K_0(\mathcal{Q}^F(\mathcal{G}))$.
	Hence they coincide by \eqref{eqn:Iwasawa-K-sequence}.
	In other words, $P$ and $\Lambda^{\mathcal{O}}(\mathcal{G})^m$ are
	stably isomorphic. By enlarging $m$ if necessary, we may assume
	that $P$ is free of rank $m$ and we have established the existence of
	\eqref{eqn:free-resolution}. By \cite[Lemma 5.3.11]{MR2392026} 
	this sequence induces an exact sequence of $\Z_p[G]$-modules
	\[
	0 \rightarrow M^{\Gamma'} \rightarrow \Z_p[G]^m \rightarrow
	\Z_p[G]^m \rightarrow M_{\Gamma'} \rightarrow 0.
	\]
	It follows that $M^{\Gamma'}$ is a free $\Z_p$-module of the same rank
	as $M_{\Gamma'}$. This proves the remaining claims.
\end{proof}

Now let $\pi: \mathcal{G} \rightarrow \Aut_{\mathcal{O}}(T_{\pi})$
be a continuous homomorphism as in \eqref{eqn:rho}
and set $n := \rank_{\mathcal{O}}(T_{\pi})$.
There is a ring homomorphism $\Phi_{\pi}: \Lambda(\mathcal{G}) \rightarrow M_{n\times n}(\Lambda^{\mathcal{O}}(\Gamma))$
induced by the continuous group homomorphism
\begin{eqnarray*}
\mathcal{G} & \longrightarrow & (M_{n \times n}(\mathcal{O}) \otimes_{\Z_p} \Lambda(\Gamma))^{\times} = \GL_{n}(\Lambda^{\mathcal{O}}(\Gamma))\\
g & \mapsto & \pi(g) \otimes \overline{g},
\end{eqnarray*}
where $\overline{g}$ denotes the image of $g$ in $\mathcal{G} / H \simeq \Gamma$. 
By \cite[Lemma 3.3]{MR2217048} this
homomorphism extends to a ring homomorphism
$\mathcal{Q}(\mathcal{G}) \rightarrow M_{n\times n}(\mathcal{Q}^{F}(\Gamma))$
and this in turn induces a homomorphism
\[
\Phi_{\pi}: K_{1}(\mathcal{Q}(\mathcal{G})) \longrightarrow 
K_{1}(M_{n\times n}(\mathcal{Q}^{F}(\Gamma))) 
\simeq \mathcal{Q}^{F}(\Gamma)^{\times}.
\]
For $\xi \in K_{1}(\mathcal{Q}(\mathcal{G}))$ we set $\xi(\pi) := \Phi_{\pi}(\xi)$.
If $\pi = \pi_{\rho}$ is an Artin representation with character $\rho$,
we also write $\Phi_{\rho}$ and $\xi(\rho)$ for $\Phi_{\pi_{\rho}}$ 
and $\xi(\pi_{\rho})$, respectively, and let
\[
	j_{\rho}: \zeta(\mathcal{Q}(\mathcal{G}))^{\times} \longrightarrow
	\mathcal{Q}^{F}(\Gamma)^{\times}
\]
be the map defined by Ritter and Weiss in \cite{MR2114937}.
By \cite[Lemma 2.3]{MR3072281} (choose $r=0$) we have a commutative triangle
\[ \xymatrix{
K_1(\mathcal{Q}(\mathcal{G})) \ar[d]^{\Nrd} \ar[dr]^{\Phi_{\rho}} & \\
\zeta(\mathcal{Q}(\mathcal{G}))^{\times} \ar[r]^{j_{\rho}} &
\mathcal{Q}^{F}(\Gamma)^{\times}
}\]
We shall also write $\xi^{\ast}(\rho)$ for the leading term at $T=0$ of the power
series $\Phi_{\rho}(\xi)$.

We choose a maximal $\Lambda(\Gamma_0)$-order $\mathfrak{M}(\mathcal{G})$
in $\mathcal{Q}(\mathcal{G})$ such that $\Lambda(\mathcal{G})$
is contained in $\mathfrak{M}(\mathcal{G})$.

\begin{lemma} \label{lem:char-elements}
Let $C^{\bullet} \in \mathcal{D}^{\perf}_{\tors}(\Lambda(\mathcal{G}))$ be a complex
and let $\xi$ be a characteristic element for $C^{\bullet}$.
Let $\pi_{\rho}$
be an Artin representation with character $\rho$.
Then $\xi(\rho) \cdot j_{\rho}(x)$ 
is a characteristic element for $C^{\bullet}_{\rho}$
for every 
$x \in \zeta(\mathfrak{M}(\mathcal{G}))^{\times}$.
\end{lemma}

\begin{proof}
In the case that $x = 1$ this is \cite[Lemma 5.4 (vii)]{MR4092926}
(and follows from the naturality of connecting homomorphisms).
By \cite[Remark H]{MR2114937} the image of $\zeta(\mathfrak{M}(\mathcal{G}))^{\times}$ under $j_{\rho}$ is contained in
$\Lambda^{\mathcal{O}}(\Gamma)^{\times}$. This proves the claim.
\end{proof}

\subsection{Cohomology with compact support} \label{subsec:compact-sc}
Let $p$ be an odd prime.
We denote the cyclotomic $\Z_p$-extension of a number field $K$ by
$K_{\infty}$ and set $\Gamma_K := \mathrm{Gal}(K_{\infty}/K)$. 
Let $S$ be a finite set of places of $K$ containing the set
$S_{\infty} \cup S_p$.
Let $M$ be a topological $G_{K,S}$-module.
Following Burns and Flach \cite{MR1884523} we define the compactly supported cohomology complex to be
\[
R\Gamma_c(\mathcal O_{K,S}, M) := \mathrm{cone} \left(R\Gamma(\mathcal O_{K,S}, M) \longrightarrow
\bigoplus_{v \in S} R\Gamma(K_v, M) \right)[-1],
\]
where the arrow is induced by the natural restriction maps. 
For any integers $i$ and $r$ we abbreviate
$H^iR\Gamma_c(\mathcal O_{K,S}, M)$ to $H^i_c(\mathcal O_{K,S}, M)$
and set 
$H^i_c(\mathcal O_{K,S}, \Q_p(r)) := \Q_p \otimes_{\Z_p} H^i_c(\mathcal O_{K,S}, \Z_p(r))$.

Let $E/K$ be a finite Galois extension of totally real fields 
and set $L := E(\zeta_p)$. Then $L$ is a CM-field
and we denote its maximal totally real subfield by $L^+$ as in \S 
\ref{sec:comparison}. Set $\mathcal{G}:= \Gal(L_{\infty}/K)$
and let
\[
\chi_{\mathrm{cyc}}: \mathcal{G} \longrightarrow \Z_{p}^{\times},
\]
be the $p$-adic cyclotomic character
defined by $\sigma(\zeta) = \zeta^{\chi_{\mathrm{cyc}}(\sigma)}$ for any $\sigma \in \mathcal{G}$ and any $p$-power root of unity $\zeta$.
The composition of $\chi_{\mathrm{cyc}}$ with the projections 
onto the first and second factors of the canonical decomposition
$\Z_{p}^{\times} = \langle \zeta_{p-1} \rangle \times (1+p\Z_{p})$
are given by the Teichm\"{u}ller character $\omega$ and a map
that we denote by $\kappa$.

Assume in addition that $S$ contains
all places that ramify in $L_{\infty}/K$.
For each integer $r$ we define a complex of
$\Lambda(\mathcal{G})$-modules
\[
	C_{r,S}^{\bullet} := R\Gamma_c(\mathcal{O}_{K,S}, e_r \Lambda(\mathcal{G})^{\sharp}(r)),
\]
where $\Lambda(\mathcal{G})^{\sharp}(r)$ denotes the $\Lambda(\mathcal{G})$-module
$\Lambda(\mathcal{G})$ upon which $\sigma \in G_K$ acts on the right via
multiplication by the element $\chi_{\mathrm{cyc}}^r(\sigma)\overline \sigma^{-1}$; here $\overline \sigma$
denotes the image of $\sigma$ in $\mathcal{G}$. Note that
the complexes $C_{r,S}^{\bullet}$ are perfect by \cite[Proposition 1.6.5]{MR2276851} and
we have natural isomorphisms
\begin{equation} \label{eqn:C_r-twists}
	C_{r,S}^{\bullet} \simeq C_{1,S}^{\bullet} \otimes_{\Z_p}^{\mathbb L} \Z_p(r-1)
\end{equation}
for every integer $r$. 

Each $\Lambda(\mathcal{G})$-module $M$ naturally
decomposes as a direct sum $M^+ \oplus M^-$ with 
$M^{\pm} = \frac{1 \pm \tau}{2} M$. Similarly, each complex $C^{\bullet}$ of 
$\mathcal{D}^{\perf}(\Lambda(\mathcal{G}))$ gives rise to subcomplexes
$(C^{\bullet})^+$ and $(C^{\bullet})^-$. 
Moreover, we let $(C^{\bullet})^{\vee} := R\Hom(C^{\bullet}, \Q_p/\Z_p)$ 
be the Pontryagin dual of $C^{\bullet}$.
By a Shapiro Lemma argument and 
Artin--Verdier duality we then have isomorphisms
\begin{equation} \label{eqn:C_r-duality}
C_{r,S}^{\bullet} \simeq \left\{ \begin{array}{lll}
R\Gamma(\mathcal{O}_{L^+_{\infty},S}, \Q_p/\Z_p(1-r))^{\vee}[-3] & \mbox{ if } & 2 \nmid r \\
\left(R\Gamma(\mathcal{O}_{L_{\infty},S}, \Q_p/\Z_p(1-r))^{\vee}\right)^-[-3]
& \mbox{ if } & 2 \mid r.
\end{array}
\right.
\end{equation}
in $\mathcal{D}(e_r\Lambda(\mathcal{G}))$.
We let $M_{S}$ be the maximal abelian pro-$p$-extension of $L_{\infty}^+$
unramified outside $S$. Then $X_{S} := \Gal(M_{S}/L_{\infty}^+)$
is a finitely generated $\Lambda(\mathcal{G}^+)$-module,
where we put $\mathcal{G}^+ := \mathcal{G}/\langle \tau \rangle = \Gal(L^+_{\infty}/K)$. Iwasawa \cite{MR0349627} has shown that
$X_{S}$ is in fact torsion as a $\Lambda(\Gamma_{L^+})$-module.
We let $\mu_{p}(L^+)$ denote the Iwasawa $\mu$-invariant 
of $X_{S}$ and note that
this does not depend on the choice of $S$ (see \cite[Corollary 11.3.6]{MR2392026}).
Hence $\mu_{p}(L^+)$ vanishes if and only if $X_{S}$ is finitely generated as a $\Z_{p}$-module.
It is conjectured that we always have $\mu_{p}(L^+)=0$ and as explained in \cite[Remark 4.3]{MR3749195}, 
this is closely related to the classical Iwasawa `$\mu=0$' conjecture
for $L$ at $p$. Thus a result of Ferrero and Washington \cite{MR528968} on this latter conjecture implies 
that $\mu_{p}(L^+)=0$ whenever $E/\Q$ and thus $L/\Q$ is abelian.

The only non-trivial cohomology groups of $R\Gamma(\mathcal{O}_{L^+_{\infty},S}, \Q_p/\Z_p)^{\vee}$ occur in degrees $-1$ and $0$ and canonically identify
with $X_{S}$ and $\Z_p$, respectively. Hence \eqref{eqn:C_r-duality}
with $r=1$ and \eqref{eqn:C_r-twists} imply that for each integer $r$
the cohomology of $C_{r,S}^{\bullet}$ is concentrated in degrees $2$ and $3$ and we have
\[
	H^2(C_{r,S}^{\bullet}) \simeq X_{S}(r-1), \qquad 
	H^3(C_{r,S}^{\bullet}) \simeq \Z_p(r-1).
\]

\subsection{The main conjecture}
The following is an obvious reformulation of the equivariant Iwasawa main
conjecture for the extension $L^+_{\infty}/K$ (without its uniqueness statement).

\begin{conj}[equivariant Iwasawa main conjecture] \label{conj:EIMC}
Let $L/K$ be a Galois CM-extension such that $L$ contains a primitive $p$-th
root of unity.
Let $S$ be a finite set of places of $K$ containing $S_{\infty}$ and all places
that ramify in $L^+_{\infty}/K$.
Then there exists an element $\zeta_{S} \in K_{1}(\mathcal{Q}(\mathcal{G}^+))$ such that
$\partial(\zeta_{S}) = [C_{1,S}^{\bullet}]$ and 
for every irreducible Artin representation $\pi_{\rho}$ of $\mathcal{G}^+$ with character
$\rho$ and for each integer $n \geq 1$ divisible by $p-1$ we have
\begin{equation}\label{eqn:EIMC-eqn}
  \zeta_{S}^{\ast}(\rho \kappa^{n}) = L_{p, S}(1-n, \rho) = j\left(L_{S}(1-n, \rho^{j^{-1}})\right) 
\end{equation}
for every field isomorphism $j: \C \cong \C_p$.
\end{conj}


Part (i) of the following theorem has been shown by 
Ritter and Weiss \cite{MR2813337} and by Kakde \cite{MR3091976} independently.
Part (ii) is due to Johnston and the present author \cite{abelian-MC}.

\begin{theorem}\label{thm:EIMC-known-cases}
Conjecture \ref{conj:EIMC} holds for $L_{\infty}^+/K$ in each of the following cases.
\begin{enumerate}
	\item 
	The $\mu$-invariant $\mu_{p}(L^+)$ vanishes.
	\item
	The Galois group $\mathcal{G}^+$ has an abelian Sylow $p$-subgroup.
\end{enumerate}
\end{theorem}

By starting out from the work of Deligne and Ribet \cite{MR579702},
Greenberg \cite{MR692344} has shown that for each topological
generator $\gamma$ of $\Gamma$ there is a unique element
$f_{\rho, S}(T)$ in the quotient field of 
$\Q_p^c \otimes_{\Z_p} \Z_p \llbracket T \rrbracket$ 
(which we will identify with $\Q_p^c \otimes_{\Z_p} \Lambda(\Gamma)$ via the usual map
that sends $\gamma$ to $1+T$)
such that
\[
	L_{p,S}(1-s, \rho) = f_{\rho, S}(u^s-1),
\]
where $u := \kappa(\gamma)$.
For each integer $a$ we let $x \mapsto t_{\cyc}^a(x)$ be the
automorphisms on $\mathcal{Q}^c(\mathcal G)$ induced by
$g \mapsto \chi_{\mathrm{cyc}}(g)^a g$ for $g \in \mathcal{G}$.
We use the same notation for the induced group homomorphisms
on $K_1(\mathcal{Q}(\mathcal{G}))$ and $K_1(e_r\mathcal{Q}(\mathcal{G}))$,
$r \in \Z$.

\begin{prop}\label{prop:rho-of-zeta}
Suppose that Conjecture \ref{conj:EIMC} holds for $L_{\infty}^+/K$.
Then for each $r \in \Z$ there exists an element 
$\zeta_{r,S} \in K_{1}(e_r\mathcal{Q}(\mathcal{G}))$ such that
$\partial(\zeta_{r,S}) = [C_{r,S}^{\bullet}]$ and 
for every irreducible Artin representation $\pi_{\rho}$ of $\mathcal{G}^+$ with character
$\rho$ we have
\begin{equation}\label{eqn:rho-of-zeta}
	\zeta_{r,S}(\rho \otimes \omega^{r-1}) = f_{\rho, S}(u^{1-r}(1+T)-1). 
\end{equation}
\end{prop}

\begin{proof}
When $r=1$ we may take $\zeta_{1,S} = \zeta_S$ by \cite[Proposition 7.5]{MR4111943}. 
Then $\zeta_{r,S} := t_{\cyc}^{1-r}(\zeta_{1,S})$ is a characteristic element for 
$C_{r,S}^{\bullet}$ by \eqref{eqn:C_r-twists} and \eqref{eqn:rho-of-zeta} follows from 
\cite[Lemma 5.4 (v)]{MR4092926} (see also \cite[Lemma 9.5]{MR3294653}).
\end{proof}

\begin{corollary} \label{cor:char-element}
Let $r \in \Z$ be an integer and let $\pi_{\rho}$ be an irreducible Artin representation  of $\mathcal{G}^+$ with character $\rho$.
Then $f_{\rho, S}(u^{1-r}(1+T)-1)$ is a characteristic
element of $C^{\bullet}_{r, \rho \otimes \omega^{r-1}}$.
\end{corollary}

\begin{proof}
Since the main conjecture
(Conjecture \ref{conj:EIMC}) holds `over the maximal order' by
\cite[Theorem 4.9]{MR3749195} (this result is essentially due
to Ritter and Weiss \cite{MR2114937}), 
the equality \eqref{eqn:rho-of-zeta} holds unconditionally up to a factor 
$j_{\rho}(x)$ for some $x \in \zeta(\mathfrak{M}(\mathcal{G}))^{\times}$.
Thus the claim follows from Lemma \ref{lem:char-elements}.
\end{proof}

\subsection{Schneider's conjecture and semisimplicity} \label{subsec:Schneider}
We recall the following result from \cite[Propositions 3.11 and
 3.12]{wild-kernels}. Part (i) is a special case of
 \cite[Proposition 1.20]{MR1386106} and of \cite[Proposition 1.6.5]{MR2276851}.

\begin{prop}  \label{prop:cs-cohomology}
 Let $L/K$ be a Galois extension of number fields with Galois group $G$.
 Let $r>1$ be an integer and let $p$ be an odd prime. Then the following hold.
 \begin{enumerate}
  \item 
  The complex $R\Gamma_c(\mathcal O_{L,S}, \Z_p(r))$
  belongs to $\mathcal D^{\perf} (\Z_p[G])$ and is acyclic outside degrees
  $1,2$ and $3$.
  \item
  We have an exact sequence of $\Z_p[G]$-modules
  \[
   0 \longrightarrow H_{-r}^+(L) \otimes \Z_p \longrightarrow H^1_c(\mathcal O_{L,S}, \Z_p(r)) \longrightarrow
   \Sha^1(\mathcal O_{L,S}, \Z_p(r)) \longrightarrow 0.
  \]
  \item
    We have an isomorphism of $\Z_p[G]$-modules
    \[
     H^3_c(\mathcal O_{L,S}, \Z_p(r)) \simeq \Z_p(r-1)_{G_L}.
    \]
    \item 
    The $\Z_{p}[G]$-module $\Sha^2(\mathcal O_{L, S}, \Z_p(r))$
    is finite.
  \item
  The $\Z_p$-rank of
  $H^2_c(\mathcal O_{L,S}, \Z_p(r))$ equals $d_{r+1}$ if and only if
  Schneider's conjecture $\Sch(L,p,r)$ holds.
 \end{enumerate}
\end{prop}

We now return to the situation considered in \S \ref{subsec:compact-sc}.
Before proving the next result we recall that for a finitely generated
$\Lambda^{\mathcal{O}}(\Gamma)$-torsion module $M$ the following
are equivalent: (i) $M^{\Gamma}$ is finite; (ii) $M_{\Gamma}$ is finite;
(iii) $f_M(0) \not= 0$, where $f_M$ denotes a characteristic element for $M$.
Moreover, by \cite[Proposition 5.3.19]{MR2392026} we have that $M=0$ if and only
if $M_{\Gamma} = M^{\Gamma} = 0$.

\begin{prop} \label{prop:descent-iso}
Let $r>1$ be an integer and let $\rho$
be an irreducible Artin character of $\mathcal{G}^+$ such that 
$\ker(\rho)$ contains $\Gamma_{L^+}$.
Set $\psi := \rho \otimes \omega^{r-1}$. Then we have
that $H^i(C^{\bullet}_{r,S,\psi}) = 0$ for $i \not= 2,3$ 
and natural isomorphisms of $\mathcal{O}_{\psi}$-modules
\begin{equation} \label{eqn:equation-1c}
	H^2(C^{\bullet}_{r,S,\psi})^{\Gamma} \simeq 
	H^1_c(\mathcal O_{L,S}, \Z_p(r))^{(\psi)}  
	\simeq \Sha^1(\mathcal{O}_{L, S}, \Z_p(r))^{(\psi)} 
\end{equation}
and
\begin{equation} \label{eqn:equation-3c}
	H^3(C^{\bullet}_{r,S,\psi})_{\Gamma} \simeq 
	H^3_c(\mathcal{O}_{L, S}, \Z_p(r))_{(\psi)}  \simeq 
	(\Z_p(r-1)_{G_L})_{(\psi)}
\end{equation}
and a short exact sequence of $\mathcal{O}_{\psi}$-modules
\begin{equation} \label{eqn:equation-2c}
	0 \longrightarrow H^2(C^{\bullet}_{r,S,\psi})_{\Gamma} \longrightarrow
	H^2(R\Gamma_c(\mathcal O_{L,S}, \Z_p(r))_{(\psi)}) \longrightarrow
	H^3(C^{\bullet}_{r,S,\psi})^{\Gamma} \longrightarrow 0.
\end{equation}
In particular, the $\mathcal{O}_{\psi}$-modules $H^3(C^{\bullet}_{r,S,\psi})_{\Gamma}$
and $H^3(C^{\bullet}_{r,S,\psi})^{\Gamma}$ are finite.
\end{prop}

\begin{proof}
We put $\mathcal{O} := \mathcal{O}_{\psi}$ for brevity.
Since the complex $C_{r,S}^{\bullet}$ is acyclic 
outside degrees $2$ and $3$ and the functor
$M \mapsto M_{\psi} = 
(\Lambda^{\mathcal{O}}(\Gamma) \otimes_{\mathcal{O}} T_{\psi}) 
\otimes_{\Lambda(\mathcal{G})} M$ is right exact, it is clear that
$H^i(C^{\bullet}_{r,S,\psi})$ vanishes for $i\geq 4$.
We let $U = \Gamma_L  =\Gal(L_{\infty}/L)$. Then $U$ is contained in
the kernel of $\psi$, and by \cite[Proposition 1.6.5]{MR2276851}
we have an isomorphism 
\[
C^{\bullet}_{r,S,U} \simeq 
e_r R\Gamma_c(\mathcal{O}_{L,S}, \Z_p(r)) :=
e_r \Z_{p}[G] \otimes_{\Z_p[G]}^{\mathbb L} R\Gamma_c(\mathcal{O}_{L,S}, \Z_p(r))
\]
in $\mathcal{D}(e_r\Z_p[G])$.
We now consider the exact sequence \eqref{eqn:tautological-ses} for
the case at hand and various integers $i$. 
We will repeatedly apply Lemma \ref{lem:chi-twists}.
In particular, the complex $C_{r,S,U,(\psi)}^{\bullet}$
is acyclic outside degrees $1$, $2$ and $3$
by Proposition \ref{prop:cs-cohomology} (i).
For $i\leq 0$ we find that
$H^{i}(C^{\bullet}_{r,S,\psi})_{\Gamma}$ and
$H^{i+1}(C^{\bullet}_{r,S,\psi})^{\Gamma}$ vanish. 
Thus $H^i(C^{\bullet}_{r,S,\psi})$
vanishes for $i \leq 0$ and even for $i=1$ once we show that
the $\mathcal{O}$-module
$H^{1}(C^{\bullet}_{r,S,\psi})_{\Gamma}$ is trivial.
We already know that it is finite.
Sequence \eqref{eqn:tautological-ses}
in the case $i=1$ and Lemma \ref{lem:chi-twists} now give 
rise to a short exact sequence
\[
0 \longrightarrow H^{1}(C^{\bullet}_{r,S,\psi})_{\Gamma} \longrightarrow
H^1_c(\mathcal O_{L,S}, \Z_p(r))^{(\psi)}  \longrightarrow
H^2(C^{\bullet}_{r,S,\psi})^{\Gamma} \longrightarrow 0.
\]
Since the central idempotent $e_r$ annihilates $H_{-r}^+(L) \otimes \Z_p$
and $e_r e_{\psi} = e_{\psi}$ we have an isomorphism
\[
H^1_c(\mathcal O_{L,S}, \Z_p(r))^{(\psi)}  
	\simeq \Sha^1(\mathcal{O}_{L, S}, \Z_p(r))^{(\psi)} 
\]
by Proposition \ref{prop:cs-cohomology} (ii).
Since the Tate--Shafarevich group $\Sha^1(\mathcal{O}_{L, S}, \Z_p(r))$
is torsion-free by Lemma \ref{lem:Sha-torsion-free}, the $\mathcal{O}$-module
$\Sha^1(\mathcal{O}_{L, S}, \Z_p(r))^{(\psi)}$ is free. Thus 
$H^{1}(C^{\bullet}_{r,S,\psi})_{\Gamma}$ vanishes as desired 
and we have established
\eqref{eqn:equation-1c}. Proposition \ref{prop:cs-cohomology} (iii),
Lemma \ref{lem:chi-twists} and
the case $i=3$ of sequence \eqref{eqn:tautological-ses} imply
\eqref{eqn:equation-3c}.
Finally, sequence \eqref{eqn:equation-2c} is the case $i=2$ of sequence
\eqref{eqn:tautological-ses}.
\end{proof}

\begin{lemma} \label{lem:explicit-complex}
Let $r$ be an arbitrary integer.
Then there are finitely generated
$e_r\Lambda(\mathcal{G})$-modules $Y_{r,S}$ and $Z_r$
with all of the following properties:
\begin{enumerate}
\item 
The projective dimension of both $Y_{r,S}$ and $Z_r$ is at most $1$;
\item 
both $Y_{r,S}$ and $Z_r$ are torsion as $R$-modules; 
\item 
there is an exact triangle
\[
	Z_r[-3] \longrightarrow C_{r,S}^{\bullet} \longrightarrow Y_{r,S}[-2]
\]
in $\mathcal{D}(e_r\Lambda(\mathcal{G}))$;
\item 
we have that $Y_{r,S} = Y_{0,S}(r)$ and $Z_r = Z_0(r)$;
\item 
the coinvariants $(Z_r)_{\Gamma_L}$ are finite if $r \not=1$.
\end{enumerate}
\end{lemma}

\begin{proof}
 We first consider the case $r=0$. It is shown in
\cite[Proposition 8.5]{BS-MathZ} that the complex $C_{0,S}^{\bullet}$
is isomorphic in $\mathcal{D}(e_0 \Lambda(\mathcal{G}))$ to a complex
\[
 \dots \longrightarrow 0 \longrightarrow Y_{0,S} \longrightarrow Z_0 
 \longrightarrow 0 \longrightarrow \dots
\]
where $Y_{0,S}$ is placed in degree $2$. More precisely, in the notation of
\cite{BS-MathZ} we have 
$C_{0,S}^{\bullet} = C_S^{\bullet}(L_{\infty}^+/K)(-1)[-3]$,
$Y_{0,S} = Y_S^T(-1)$ and
$Z_0  = I_T = \left(\bigoplus_{v \in T} \ind_{\mathcal{G}_{w_{\infty}}}^{\mathcal{G}}
\Z_p(-1)\right)^-$. Here $T$ is a finite set of places of $K$ disjoint from $S$
with certain properties, and $\mathcal{G}_{w_{\infty}}$ denotes the decomposition
group at a chosen place $w_{\infty}$ of $L_{\infty}$ above $v$
for each $v \in T$; moreover,
we write $\ind_{\mathcal{U}}^{\mathcal{G}} M := \Lambda(\mathcal{G}) 
\otimes_{\Lambda(\mathcal{U})} M$ for any open subgroup $\mathcal{U}$ 
of $\mathcal{G}$ and any $\Lambda(\mathcal{U})$-module $M$.
By \cite[Lemmas 8.4 and 8.5]{BS-MathZ} the modules $Z_0$ and $Y_{0,S}$
are $R$-torsion and of projective dimension at most $1$. 
Thus (i) and (ii) also hold
for $Y_{r,S} := Y_{0,S}(r) = Y_S^T(r-1)$ and $Z_r := Z_0(r) = e_r
\bigoplus_{v \in T} \ind_{\mathcal{G}_{w_{\infty}}}^{\mathcal{G}}
\Z_p(r-1)$. It is now clear that (iv) holds and that
$C_{r,S}^{\bullet}$ is isomorphic to the complex
\[
 \dots \longrightarrow 0 \longrightarrow Y_{r,S} \longrightarrow Z_r
 \longrightarrow 0 \longrightarrow \dots
\]
in $\mathcal{D}(e_r \Lambda(\mathcal{G}))$. Hence (iii) also holds.
Moreover, the coinvariants $(Z_r)_{\Gamma_L}$ are clearly finite for $r \not=1$.
\end{proof}

\begin{remark} \label{rem:similar-constructions}
We point out that similar constructions repeatedly appear in the literature.
In fact, by \cite[Theorem 2.4]{MR3072281} the complex $C_S^{\bullet}(L_{\infty}^+/K)$ naturally identifies with the complex constructed
by Ritter and Weiss \cite{MR1935024}. Choosing their maps $\Psi$ and $\tilde\psi$
in \cite[p.562 f.]{MR2114937} suitably one 
can take $Y_{1,S} = \coker(\Psi)$ and $Z_1 = \coker(\tilde \psi)$.
Moreover, Burns \cite[\S 5.3.1]{MR4092926} constructed an exact triangle
in $\mathcal{D}(e_0 \Lambda(\mathcal{G}))$ of the form
\[
	R\Gamma_T(\mathcal{O}_{K,S}, e_0 \Lambda(\mathcal{G})^{\sharp}(1)) 
	\longrightarrow R\Gamma(\mathcal{O}_{K,S}, e_0 \Lambda(\mathcal{G})^{\sharp}(1)) 
	\longrightarrow \bigoplus_{v \in T} R\Gamma(K(v), e_0 \Lambda(\mathcal{G})^{\sharp}(1)),
\]
where the set $T$ is as in the proof of Lemma \ref{lem:explicit-complex}.
The complexes $R\Gamma(K(v), e_0 \Lambda(\mathcal{G})^{\sharp}(1))$ are acyclic outside
degree $1$ and we have $e_0 \Lambda(\mathcal{G})$-isomorphisms
\[
	H^1(K(v), e_0 \Lambda(\mathcal{G})^{\sharp}(1)) \simeq 
	e_0 \ind_{\mathcal{G}_{w_{\infty}}}^{\mathcal{G}} \Z_p(1)	
\]
for each $v \in T$.
For each $\Lambda(\mathcal{G})$-module $M$ and $i \in \Z$ we set
$E^i(M) := \Ext^i_{\Lambda(\mathcal{G})}(M, \Lambda(\mathcal{G}))$. We then
have an isomorphism of $e_0 \Lambda(\mathcal{G})$-modules
\[
	\bigoplus_{v \in T} H^1(K(v), e_0 \Lambda(\mathcal{G})^{\sharp}(1)) \simeq 
	E^1(I_T).
\]
 If $C^{\bullet}$ is a complex in $\mathcal{D}(e_0\Lambda(\mathcal{G}))$,
we write $(C^{\bullet})^{\ast}$ for the complex 
$R\Hom_{e_0\Lambda(\mathcal{G})}(C^{\bullet}, e_0 \Lambda(\mathcal{G}))$.
If $M$ is an $R$-torsion $e_0\Lambda(\mathcal{G})$-module 
of projective dimension at most $1$,
then we have isomorphisms $M[-n]^{\ast} \simeq E^1(M)[n-1]$ in
$\mathcal{D}(e_0\Lambda(\mathcal{G}))$ for every $n \in \Z$. This yields
\[
\bigoplus_{v \in T} R\Gamma(K(v), e_0 \Lambda(\mathcal{G})^{\sharp}(1))
\simeq I_T^{\ast}.
\]
The complex $R\Gamma_T(\mathcal{O}_{K,S}, e_0 \Lambda(\mathcal{G})^{\sharp}(1))$
is acyclic outside degree $2$ and
its second cohomology group
$H^2_T(\mathcal{O}_{K,S}, e_0 \Lambda(\mathcal{G})^{\sharp}(1))$
is of projective dimension at most $1$ by \cite[Proposition 5.5]{MR4092926}.
Since we have an isomorphism
\[
	R\Gamma(\mathcal{O}_{K,S}, e_0 \Lambda(\mathcal{G})^{\sharp}(1)) \simeq
	(C_{0,S}^{\bullet})^{\ast}[-3]
\]
by Artin--Verdier duality, it follows that we have 
an isomorphism
\[
R\Gamma_T(\mathcal{O}_{K,S}, e_0 \Lambda(\mathcal{G})^{\sharp}(1)) \simeq 
Y_S^T(-1)^{\ast}[-1]
\]
in $\mathcal{D}(e_0\Lambda(\mathcal{G}))$ and an isomorphism of
$e_0 \Lambda(\mathcal{G})$-modules
\[
	H^2_T(\mathcal{O}_{K,S}, e_0 \Lambda(\mathcal{G})^{\sharp}(1)) \simeq 
	E^1(Y_S^T(-1)).
\]
Finally, the minus Tate module $T_p(\mathcal{M_{S,T}^K})^-$ 
of the Iwasawa-theoretic $1$-motive $\mathcal{M_{S,T}^K}$
constructed by Greither and Popescu \cite{MR3383600} may play the role
of $E^1(Y_S^T(-1))$; the $e_0 \Lambda(\mathcal{G})$-module
$T_p(\Delta_{\mathcal{K,T}})^-$ is our $E^1(I_T)$. See in particular
\cite[Remark 3.10]{MR3383600}.
\end{remark}

\begin{theorem} \label{thm:Schneider-semisimple}
Let $r>1$ be an integer and let $\rho$
be an irreducible Artin character of $\mathcal{G}^+$ such that 
$\ker(\rho)$ contains $\Gamma_{L^+}$.
Set $\psi := \rho \otimes \omega^{r-1}$. Then the following
conditions are equivalent.
\begin{enumerate}
\item 
We have that $\Sha^1(\mathcal{O}_{L, S}, \Z_p(r))^{(\psi)}$ vanishes.
\item 
We have that $\Sha^1(\mathcal{O}_{L, S}, \Q_p(r))^{(\psi)}$ vanishes.
\item 
We have that $H^2_c(\mathcal{O}_{L, S}, \Q_p(r))^{(\psi)}$ vanishes.
\item
The period $\Omega_j (r, \psi)$ is non-zero for any (and hence every) choice
of $j: \C \simeq \C_p$.
\item 
We have that $L_{p, S}(r, \rho) \not= 0$.
\end{enumerate}
If these equivalent conditions hold (in particular if $\Sch(L,p,r)$ holds),
then the complex $C_{r,S}^{\bullet}$ is semisimple at $\psi$.
\end{theorem}

\begin{proof}
We have already observed in the proof of Proposition \ref{prop:descent-iso}
that $\Sha^1(\mathcal{O}_{L, S}, \Z_p(r))^{(\psi)}$ is a free 
$\mathcal{O}$-module, where we again set $\mathcal{O} := \mathcal{O}_{\psi}$.
The equivalence of (i) and (ii) is therefore clear.
We have an exact sequence of $\Z_{p}[G]$-modules
\[ \begin{array}{rllllll}
    0 & \longrightarrow & \Sha^1(\mathcal O_{L, S}, \Z_p(r)) & \longrightarrow & H^1_{\et}(\mathcal O_{L, S}, \Z_p(r)) 
    & \longrightarrow & P^1(\mathcal O_{L, S}, \Z_p(r))\\
    & \longrightarrow & H^2_c(\mathcal O_{L, S}, \Z_p(r)) & \longrightarrow & \Sha^2(\mathcal O_{L, S}, \Z_p(r)) & \longrightarrow & 0.
    \end{array}
\]
Since the $e_r\Q_p[G]$-modules
$e_r H^1_{\et}(\mathcal O_{L, S}, \Q_p(r)) \simeq
e_r K_{2r-1}(\mathcal{O}_L) \otimes \Q_p$ and $e_r P^1(\mathcal O_{L, S}, \Q_p(r))
\simeq e_r(L \otimes_{\Q} \Q_p)$ are (non-canonically) isomorphic and
$\Sha^2(\mathcal O_{L, S}, \Z_p(r))$ is finite by
Proposition \ref{prop:cs-cohomology} (iv), there is
a (non-canonical) isomorphism of $e_r \Q_p[G]$-modules
\[
	e_r \Sha^1(\mathcal{O}_{L, S}, \Q_p(r)) \simeq 
	e_r H^2_c(\mathcal{O}_{L, S}, \Q_p(r)).
\]
Thus also $\Sha^1(\mathcal{O}_{L, S}, \Q_p(r))^{(\psi)}$ and
$H^2_c(\mathcal{O}_{L, S}, \Q_p(r))^{(\psi)}$ are (non-canonically)
isomorphic and so (ii) and (iii) are indeed equivalent.
The equivalence of (ii) and (iv) is easy (see Remark \ref{rmk:rho-parts-of-Leo}).

We next establish the equivalence of (i) and (v). 
By Lemma \ref{lem:one-degree-pi} the triangle of Lemma \ref{lem:explicit-complex} (iii)
induces an exact triangle
\begin{equation} \label{eqn:mc-triangle-psi}
	Z_{r,\psi}[-3] \longrightarrow C^{\bullet}_{r,S,\psi} \longrightarrow Y_{r,S,\psi}[-2]
\end{equation}
in $\mathcal{D}(\Lambda^{\mathcal{O}}(\Gamma))$.
Let $h_{r, \psi}(T)$ and $g_{r,S, \psi}(T)$ be
characteristic elements of $Z_{r,\psi}$ and $Y_{r,S,\psi}$, respectively.
Corollary \ref{cor:char-element} implies that we may assume that
\[ 
	g_{r,S,\psi}(T) = h_{r,\psi}(T) \cdot f_{\rho,S}(u^{1-r}(1+T)-1).
\]
Since $(Z_{r,\psi})_{\Gamma}$ is finite by Lemma \ref{lem:explicit-complex} (v) 
we have that $h_{r,\psi}(0) \not= 0$.
Thus $L_{p,S}(r, \rho) = f_{\rho,S}(u^{1-r}-1)$ is non-zero if and only if
$g_{r,S,\psi}(0) \not= 0$ if and only if $Y_{r,S,\psi}^{\Gamma}$ vanishes,
where the latter equivalence uses Lemma \ref{lem:free-resolution}
(with $M = Y_{r,S,\psi}$ and $\mathcal{G} = \Gamma$)
and Lemma \ref{lem:explicit-complex} (i) and (ii).
By \eqref{eqn:mc-triangle-psi} we have an exact sequence of
$\Lambda^{\mathcal{O}}(\Gamma)$-modules
\[
0 \longrightarrow H^2(C^{\bullet}_{r,S,\psi}) \longrightarrow Y_{r,S,\psi}
\longrightarrow Z_{r,\psi} \longrightarrow H^3(C^{\bullet}_{r,S,\psi})
\longrightarrow 0.
\]
Since taking $\Gamma$-invariants is left exact and $Z_{r,\psi}^{\Gamma}$
vanishes by another application of Lemma \ref{lem:free-resolution}, 
it follows that there is an isomorphism $Y_{r,S,\psi}^{\Gamma} \simeq H^2(C^{\bullet}_{r,S,\psi})^{\Gamma}$.
The latter identifies with 
$\Sha^1(\mathcal{O}_{L, S}, \Z_p(r))^{(\psi)} $ by
Proposition \ref{prop:descent-iso} which proves the claim.

Finally, if these equivalent conditions all hold, then Proposition
\ref{prop:descent-iso} implies that the $\mathcal{O}$-modules
$H^i(C^{\bullet}_{r,S,\psi})^{\Gamma}$ and $H^i(C^{\bullet}_{r,S,\psi})_{\Gamma}$
are finite for all $i \in \Z$. 
It follows from the short exact sequences \eqref{eqn:tautological-ses}
that $H^i(C^{\bullet}_{r,S,U})$ is finite for all $i \in \Z$,
where we set $U = \Gamma_L$ as before.
Thus the whole complex
$\Q_p \otimes_{\Z_{p}} \Delta_{C_{r,S}^{\bullet},\psi}$ vanishes.
In particular, the complex $C_{r,S}^{\bullet}$ is semisimple at $\psi$.
\end{proof}

\begin{remark} \label{rem:Y-coinvariants}
	By specializing $K = L^+$ in the above argument, we see that
	$\Sch(L,p,r)$ implies that $Y_{r,S}^{\Gamma_L}$ vanishes and that 
	$(Y_{r,S})_{\Gamma_L}$ is finite.
\end{remark}

\subsection{Higher refined $p$-adic class number formulae} \label{subsec:higher-CNF}

We keep the notation of \S \ref{subsec:compact-sc}.
We let $\partial_p: \zeta(\C_p[G])^{\times} \simeq K_1(\C_{p}[G]) \rightarrow
K_0(\Z_p[G], \C_p)$ be the composition of the inverse of the reduced norm and
the connecting homomorphism to relative $K$-theory. By abuse of notation we shall
use the same symbol for the induced maps on `$e_r$-parts'.
We recall that $G$ denotes $\Gal(L/K)$ and that $G^+ = G / \langle \tau \rangle
= \Gal(L^+/K)$. We define
\[
	L_{p,S}(r) := \sum_{\rho \in \Irr_{\C_{p}}(G^+)} L_{p,S}(r, \rho) 
	e_{\rho \otimes \omega^{r-1}}
	\in \zeta(e_r\Q_p[G]).
\]
Let us assume that Schneider's conjecture $\Sch(L,p,r)$ holds.
Then we actually have that $L_{p,S}(r) \in \zeta(e_r\Q_p[G])^{\times}$ and
the cohomology groups of the complex
$e_r R\Gamma_c(\mathcal{O}_{L,S}, \Z_p(r))$
are finite by Proposition \ref{prop:cs-cohomology} 
and Theorem \ref{thm:Schneider-semisimple}.
This complex therefore is an object in $\mathcal{D}^{\perf}_{\tors}(e_r \Z_p[G])$.
We now state our conjectural higher refined $p$-adic class number formula.

\begin{conj} \label{conj:higher-p-adic-formula}
Let $r>1$ be an integer and assume that $\Sch(L,p,r)$ holds. Then in
$K_0(e_r\Z_p[G], \Q_p)$ one has
\[
	\partial_p(L_{p,S}(r)) = [e_r R\Gamma_c(\mathcal{O}_{L,S}, \Z_p(r))].
\]
\end{conj}

\begin{remark}
In the case $r=1$ Burns \cite[Conjecture 3.5]{MR4092926} has formulated
a conjectural refined $p$-adic class number formula. Conjecture 
\ref{conj:higher-p-adic-formula} might be seen as a higher analogue of his conjecture.
Accordingly, Theorem \ref{thm:p-class-holds} below is the higher analogue of
\cite[Theorem 3.6]{MR4092926}.
\end{remark}

\begin{lemma}
Let $r>1$ be an integer and assume that $\Sch(L,p,r)$ holds.
Then 
\[
	HR_p(L/K,r) := \partial_p(L_{p,S}(r)) - [e_r R\Gamma_c(\mathcal{O}_{L,S}, \Z_p(r))]
\]
does not depend on the set $S$.
\end{lemma}

\begin{proof}
Let $S'$ be a second sufficiently large finite set of places of $K$. By embedding
$S$ and $S'$ into the union $S \cup S'$ we may and do assume that $S \subseteq S'$.
By induction we may additionally assume that $S' = S \cup \left\{v \right\}$,
where $v$ is not in $S$. In particular, $v$ is unramified in $L/K$ and 
$v \nmid p$.  By \cite[(30)]{MR1884523} we have an exact triangle
\begin{equation} \label{eqn:f-triangle}
  \bigoplus_{w \mid v} R\Gamma_f(L_w, \Z_p(r))[-1] \longrightarrow R\Gamma_c(\mathcal O_{L,S'}, \Z_p(r))
 \longrightarrow R\Gamma_c(\mathcal O_{L,S}, \Z_p(r)),
\end{equation}
where $R\Gamma_f(L_w, \Z_p(r))$ is a perfect 
complex of $\Z_p[G_w]$-modules which is naturally quasi-isomorphic to
\begin{equation} \label{eqn:f-quasi-iso}
\xymatrix{
 \Z_p[G_w] \ar[rr]^{1 - \phi_w N(v)^{-r}} & & \Z_p[G_w]
}
\end{equation}
with terms in degree $0$ and $1$. We set
\[
 \varepsilon_v(r) := \left( \det\nolimits_{\C_{p}}(1 - \phi_w N(v)^{-r} \mid V_{\chi})^{-1}
 \right)_{\chi \in \Irr_{\C_p}(G)} \in \zeta(\Q_p[G])^{\times}.
\]
We compute
\begin{eqnarray*}
[e_r R\Gamma_c(\mathcal{O}_{L,S}, \Z_p(r))] - [e_r R\Gamma_c(\mathcal{O}_{L,S'}, \Z_p(r))]
& = & [e_r\bigoplus_{w \mid v} R\Gamma_f(L_w, \Z_p(r))] \\
& = & \partial_p(e_r \varepsilon_v(r))\\
& = & \partial_p(L_{p,S}(r)) - \partial_p(L_{p,S'}(r)),
\end{eqnarray*}
where the first and second equality follow from \eqref{eqn:f-triangle} and
\eqref{eqn:f-quasi-iso}, respectively. 
This implies the claim.
\end{proof}

Our main evidence for Conjecture \ref{conj:higher-p-adic-formula} 
is provided by the following result which, crucially,
does not depend upon the vanishing of $\mu_p(L^+)$.

\begin{theorem} \label{thm:p-class-holds}
Let $r>1$ be an integer and assume that $\Sch(L,p,r)$ holds.
If the equivariant Iwasawa main conjecture (Conjecture \ref{conj:EIMC})
holds for $L^+_{\infty}/K$ (and so in particular if $\mu_p(L^+) = 0$
or if $\mathcal{G}^+$ has an abelian Sylow $p$-subgroup), then
Conjecture \ref{conj:higher-p-adic-formula} holds.
\end{theorem}

\begin{proof}
We first observe that the complex $C_{r,S}^{\bullet}$ is semisimple at
$\rho \otimes \omega^{r-1}$ for all $\rho \in \Irr_{\C_{p}}(G^+)$ by
Theorem \ref{thm:Schneider-semisimple}. 
Moreover, if we put $U = \Gamma_L$ as above, then we have an isomorphism
\begin{equation} \label{eqn:complex-descents}
C^{\bullet}_{r,S,U} \simeq e_r R\Gamma_c(\mathcal O_{L,S}, \Z_p(r))
\end{equation}
in $\mathcal{D}(e_r \Z_p[G])$.
If Conjecture \ref{conj:EIMC} holds,
then by Proposition \ref{prop:rho-of-zeta}
there is a characteristic element $\zeta_{r,S}$ of $[C_{r,S}^{\bullet}]$
such that $\zeta_{r,S}(\rho\otimes \omega^{r-1}) = f_{\rho,S}(u^{1-r}(1+T)-1)$.
Since $f_{\rho,S}(u^{1-r}-1) = L_{p,S}(r, \rho) \not=0$ we have that
$\zeta_{r,S}^{\ast}(\rho\otimes \omega^{r-1}) = L_{p,S}(r, \rho)$.
If $\mu_p(L^+)$ vanishes or $p$ does not divide the cardinality of $G^+$,
then \cite[Theorem 2.2]{MR2749572} implies the claim (as noted above
the complexes 
$\Q_p \otimes_{\Z_{p}} \Delta_{C_{r,S}^{\bullet},\rho \otimes \omega^{r-1}}$
all vanish).\\
In order to avoid these assumptions, we proceed as follows.
Observe that both $(Z_r)_{\Gamma_L}$ and $(Y_{r,s})_{\Gamma_L}$ are finite
by Lemma \ref{lem:explicit-complex} (v)
and Theorem \ref{thm:Schneider-semisimple}
(or rather Remark \ref{rem:Y-coinvariants}), respectively,
since Schneider's conjecture holds by assumption. Recall from 
Lemma \ref{lem:explicit-complex} (iii) that we have an exact triangle
\[
Z_r[-3] \longrightarrow C_{r,S}^{\bullet} \longrightarrow Y_{r,S}[-2]
\]
in $\mathcal{D}(e_r\Lambda(\mathcal{G}))$. It now follows from
\eqref{eqn:complex-descents} and Lemma \ref{lem:free-resolution}
for both $Z_r$ and $Y_{r,S}$ that we likewise have an exact triangle
\[
(Z_r)_{\Gamma_L}[-3] \longrightarrow e_r R\Gamma_c(\mathcal O_{L,S}, \Z_p(r)) 
\longrightarrow (Y_{r,S})_{\Gamma_L}[-2]
\]
in $\mathcal{D}(e_r\Z_p[G])$. Let $H_r$ and $G_{r,S}$ in
$\zeta(e_r\mathcal{Q}(\mathcal{G}))^{\times}$ be the reduced norms of
characteristic elements of $Z_r$ and $Y_{r,S}$, respectively.
Note that both $H_r$ and $G_{r,S}$ are actually reduced norms of
matrices with coefficients in $e_r \Lambda(\mathcal{G})$.
Since Conjecture \ref{conj:EIMC} holds by assumption, we may assume that
$\Nrd(\zeta_{r,S}) = G_{r,S} / H_r$, where $\zeta_{r,S}$ is the characteristic element
of $[C_{r,S}^{\bullet}]$ that occurs in Proposition \ref{prop:rho-of-zeta}.
Now by (the proof of) \cite[Theorem 6.4]{MR2609173} one has
\[
	\partial_p(\overline{G_{r,S}}) = [(Y_{r,S})_{\Gamma_L}], \quad
	\partial_p(\overline{H_r}) = [(Z_{r})_{\Gamma_L}],
\]
where
\[
	\overline{G_{r,S}} := \sum_{\rho \in \Irr_{\C_{p}}(G^+)} 
	\mathrm{aug}_{\Gamma}(j_{\rho \otimes \omega^{r-1}}(G_{r,S})) e_{\rho \otimes \omega^{r-1}} \in \zeta(e_r\Q_p[G])^{\times}
\]
and $\overline{H_r}$ is defined similarly. Hence we obtain 
\begin{eqnarray*}
	[e_r R\Gamma_c(\mathcal O_{L,S}, \Z_p(r))] & = & 
		[(Y_{r,S})_{\Gamma_L}] - [(Z_{r})_{\Gamma_L}] \\
		& = & \partial_p(\overline{G_{r,S}} / \overline{H_r})\\
		& = & \partial_p\left(\sum_{\rho \in \Irr_{\C_{p}}(G^+)} 
		\mathrm{aug}_{\Gamma}(j_{\rho \otimes \omega^{r-1}}(\Nrd(\zeta_{r,S}))) e_{\rho \otimes \omega^{r-1}}\right)\\
		& \stackrel{(\ast)}{=}& \partial_p(L_{p,S}(r)).
\end{eqnarray*}
It remains to justify the last equality $(\ast)$. For this we compute
\begin{eqnarray*}
	\mathrm{aug}_{\Gamma}(j_{\rho \otimes \omega^{r-1}}(\Nrd(\zeta_{r,S}))) & = &
	\mathrm{aug}_{\Gamma}(\zeta_{r,S}(\rho \otimes \omega^{r-1}))\\
	& = & f_{\rho,S}(u^{1-r}-1)\\
	& = & L_{p,S}(r, \rho).
\end{eqnarray*}
Here the first and second equality follow from \cite[Lemma 2.3]{MR2822866}
and \eqref{eqn:rho-of-zeta}, respectively. This finishes the proof of $(\ast)$.
\end{proof}

Let us write $K_0(e_r\Z_p[G], \Q_p)_{\tors}$ for the torsion subgroup of
$K_0(e_r\Z_p[G], \Q_p)$.

\begin{corollary} \label{cor:torsion}
Let $r>1$ be an integer and assume that $\Sch(L,p,r)$ holds. 
Then we have that
$HR_p(L/K,r) \in K_0(e_r\Z_p[G], \Q_p)_{\tors}$.
\end{corollary}

\begin{proof}
Let $K' \subseteq L^+$ be a totally real field containing $K$. 
Denote the Galois group $\Gal(L/K')$ by $G'$. 
Then $HR_p(L/K,r)$ maps to $HR_p(L/K',r)$
under the natural restriction map
$K_0(e_r\Z_p[G], \Q_p) \rightarrow K_0(e_r\Z_p[G'], \Q_p)$.
If $L' \subseteq L$ is a Galois CM-extension of $K$, then likewise
$HR_p(L/K,r)$ maps to $HR_p(L'/K,r)$ under the natural quotient map
$K_0(e_r\Z_p[G], \Q_p) \rightarrow K_0(e_r\Z_p[\Gal(L'/K)], \Q_p)$.
Since $HR(L'/K',r)$ vanishes for each intermediate Galois CM-extension
of degree prime to $p$ by Theorem \ref{thm:p-class-holds}, we deduce the result
by a slight modification of the argument given in \cite[Proof of Proposition 11]{MR1423032}
(also see \cite[Proof of Corollary 2]{MR2805422}).
\end{proof}

\subsection{An application to the equivariant Tamagawa number conjecture} \label{subsec:application}
Let $r$ be an integer and let $L/K$ be a finite Galois extension
of number fields with Galois group $G$. 
We set $\Q(r)_L := h^0(\Spec(L))(r)$ which we regard
as a motive defined over $K$ and with coefficients in the semisimple
algebra $\Q[G]$.
The ETNC \cite[Conjecture 4 (iv)]{MR1884523} for the pair
$(\Q(r)_L, \Z[G])$ asserts that a certain canonical element
$T\Omega(\Q(r)_L, \Z[G])$ in $K_0(\Z[G], \R)$ vanishes.
Note that in this case the element $T\Omega(\Q(r)_L, \Z[G])$
is indeed well-defined as observed in \cite[\S 1]{MR1981031}.
If $T\Omega(\Q(r)_L, \Z[G])$ is rational, i.e. belongs to $K_0(\Z[G], \Q)$, then by means of
the canonical isomorphism
\[
	K_0(\Z[G], \Q) \simeq \bigoplus_{p} K_0(\Z_p[G], \Q_p)
\]
we obtain elements 
$T\Omega(\Q(r)_L, \Z[G])_p$ in $K_0(\Z_p[G], \Q_p)$.

If $r>1$ is an integer and $L/K$ is a Galois CM-extension,
the following result provides a strategy for proving the ETNC for
the pair $(\Q(r)_L, e_r \Z[\frac{1}{2}][G])$. We let
$T\Omega(\Q(r)_L, e_r\Z[\half][G])$ be the image of $T\Omega(\Q(r)_L, \Z[G])$
under the canonical maps
\[
K_0(\Z[G], \R) \longrightarrow K_0(\Z[\half][G], \R)
\longrightarrow K_0(e_r \Z[\half][G], \R)
\]
induced by extension of scalars.
We define $T\Omega(\Q(r)_L, e_r\Z[\half][G])_p$ similarly.

\begin{theorem} \label{thm:ETNC-at-r}
Let $r>1$ be an integer and let $p$ be an odd prime.
Let $L/K$ be a Galois CM-extension with Galois group $G$ and set
$\tilde L := L(\zeta_p)$. Assume that both
Schneider's conjecture $\Sch(\tilde L,p,r)$ and the $p$-adic Beilinson conjecture
(Conjecture \ref{conj:p-adic-Beilinson}) for $\tilde L^+/K$ hold.
Then $T\Omega(\Q(r)_L, e_r\Z[\half][G])$ is rational and we have that
\[
	T\Omega(\Q(r)_L, e_r\Z[\half][G])_p \in K_0(e_r \Z_p[G], \Q_p)_{\tors}.
\]
If we assume in addition that the equivariant Iwasawa main conjecture
(Conjecture \ref{conj:EIMC}) holds for $\tilde L^+_{\infty}/K$, 
then the $p$-part of the ETNC for
the pair $(\Q(r)_L, e_r \Z[\frac{1}{2}][G])$ holds,
i.e.~the element $T\Omega(\Q(r)_L, e_r \Z[\half][G])_p$ vanishes.
\end{theorem}

\begin{proof}
Since $T\Omega(\Q(r)_{\tilde L}, e_r\Z[\half][\Gal(\tilde L / K)])$ maps to
$T\Omega(\Q(r)_L, e_r\Z[\half][G])$ under the canonical quotient map
$K_0(e_r \Z[\half][\Gal(\tilde L / K)], \R) \rightarrow K_0(e_r \Z[\half][G], \R)$
by \cite[Theorem 4.1]{MR1884523}, we may and do assume that $\tilde L = L$.
Since both Schneider's conjecture and the $p$-adic Beilinson conjecture hold,
Corollary \ref{cor:p-Beilinson-implies-Gross} implies that
Gross's conjecture at $s=1-r$ holds for all even (odd) irreducible characters of $G$
if $r$ is even (odd). By \cite[Proposition 5.5 and Theorem 6.5]{wild-kernels}
(or rather the `$e_r$-parts' of these results) this is indeed equivalent
to the rationality of $T\Omega(\Q(r)_L, e_r\Z[\half][G])$.

Let us define
\[
	\Omega_j(r) := \sum_{\rho \in \Irr_{\C_p}(G^+)} \Omega_j(r, \rho\otimes \omega^{r-1}) 
	e_{\rho \otimes \omega^{r-1}}
	\in \zeta(e_r\C_p[G])^{\times}.
\]
By the validity of the $p$-adic Beilinson conjecture we have that
\begin{equation} \label{eqn:Beilinson-cons}
L_{p,S}(r) = j(e_r L_S^{\ast}(r)) \Omega_j(r).
\end{equation}
We clearly have that $\Omega_j(r) = \Nrd ([e_r(L \otimes_{\Q} \C_p) \mid \mu_p(r) \circ
(\C_p \otimes_{\C,j} \mu_{\infty}(r))^{-1}])$.
Moreover, by Proposition \ref{prop:crucial} the automorphism
\[
t(r,S,j) := \iota_r \circ (\C_p \otimes_{\C,j} \mu_{\infty}(r)) \circ \mu_p(r)^{-1}
\circ \iota_r^{-1} \in \Aut_{e_r \C_p[G]}(H_{1-r}^+ \otimes_{\Z} \C_p)
\]
coincides with the `$e_r$-part' of the trivialization of the same name in
\cite[\S 6.2]{wild-kernels} (up to an insignificant factor $2$; 
cf.~Remark \ref{rem:regulators}). 
Since we have that
$\Nrd([H_{1-r}^+ \otimes_{\Z} \C_p \mid t(r,S,j)]) = \Omega_{j}(r)^{-1}$, the object
\begin{equation} \label{eqn:Psi}
	\Psi_{r,S}^j := [e_r R\Gamma_c(\mathcal{O}_{L,S}, \Z_p(r))] - 
	\partial_p(\Omega_{j}(r)) \in K_0(e_r \Z_p[G], \C_p)
\end{equation}
is equal to the $e_r$-part of the object denoted by $\Omega_{r,S}^j$ in
\cite[\S 6.2]{wild-kernels}. We now compute
\begin{eqnarray*}
HR_p(L/K,r) & = & \partial_p(L_{p,S}(r)) - [e_r R\Gamma_c(\mathcal{O}_{L,S}, \Z_p(r))] \\
  & = & \partial_p(j(e_r L_{S}^{\ast}(r))) + \partial_p(\Omega_{j}(r)) 
  	- [e_r R\Gamma_c(\mathcal{O}_{L,S}, \Z_p(r))]\\
  & = & \partial_p(j(e_r L^{\ast}_{S}(r))) - \Psi_{r,S}^j\\
  & = & T\Omega(\Q(r)_L, e_r \Z[\half][G])_p.
\end{eqnarray*}
Here, the first equality holds by definition of $HR_p(L/K,r)$, 
the second and third equality follow from
\eqref{eqn:Beilinson-cons} (essentially the $p$-adic Beilinson conjecture)
and \eqref{eqn:Psi}, respectively, and the last equality follows from
\cite[Proposition 6.4 and Theorem 6.5]{wild-kernels}. 
Now Corollary \ref{cor:torsion} and Theorem \ref{thm:p-class-holds}
give the result.
\end{proof}

\begin{remark} \label{rem:lower-L}
Fix an integer $r>1$.
We have assumed throughout that $\tilde L$ contains a $p$-th root of unity.
What was actually needed in the above considerations, however, is that
$\omega^{r-1}$ restricted to $G_{\tilde L}$ is trivial. 
Let $E \subseteq \tilde L_r \subseteq \tilde L$
be the smallest intermediate field such that
$\omega^{r-1}$ restricted to $G_{\tilde L_r}$ is trivial. 
Then we can replace $\tilde L$ by
$\tilde L_r$ throughout. 
Note that $\tilde L_r$ is totally real and thus $\tilde L_r^+ = \tilde L_r$
whenever $r$ is odd, whereas 
we have $\tilde L_r = \tilde L$ otherwise.
In particular, we can replace $\tilde L$ by $E$ whenever 
$r \equiv 1 \mod (p-1)$.
\end{remark}

\begin{corollary} \label{cor:ETNC-at-r}
Let $p$ be an odd prime and let $r>1$ be an integer such that $r \equiv 1 \mod (p-1)$.
Let $E/K$ be a Galois extension of totally real fields with Galois group $G$.
Assume that Schneider's conjecture $\Sch(E,p,r)$, the $p$-adic Beilinson conjecture
(Conjecture \ref{conj:p-adic-Beilinson}) for $E/K$
and the equivariant Iwasawa main conjecture
(Conjecture \ref{conj:EIMC}) for $E_{\infty}/K$ all hold.
Then $T\Omega(\Q(r)_E, \Z[\half][G])$ is rational and the $p$-part of the ETNC for
the pair $(\Q(r)_E, \Z[\frac{1}{2}][G])$ holds.
\end{corollary}

\begin{proof}
This follows from Theorem \ref{thm:ETNC-at-r} and Remark \ref{rem:lower-L}.
\end{proof}

\begin{corollary} \label{cor:classical-result}
Let $L/K$ be a Galois extension of number fields with Galois group $G$ and 
let $p$ be an odd prime. Assume that $L/\Q$
is abelian. Then $T\Omega(\Q(r)_L, e_r\Z[\half][G])$ is rational 
and the $p$-part of the ETNC for
the pair $(\Q(r)_L, e_r\Z[\frac{1}{2}][G])$ holds for all but finitely many $r>1$.
\end{corollary}

\begin{proof}
We may assume that $K = \Q$ by functoriality.
As $L(\zeta_p)$ is abelian over the rationals, the relevant 
Iwasawa invariant $\mu_p(L(\zeta_p)^+)$
vanishes by the aforementioned result of Ferrero and Washington \cite{MR528968}
(see the discussion following \eqref{eqn:C_r-duality}). Hence Conjecture
\ref{conj:EIMC} holds for $L(\zeta_p)^+_{\infty}/\Q$ by either part of Theorem \ref{thm:EIMC-known-cases} 
(but note that in this case a variant of the equivariant Iwasawa main conjecture 
can be deduced from
work of Mazur and Wiles \cite{MR742853} as in \cite[Theorem 8]{MR1935024}
for example). The $p$-adic Beilinson conjecture holds
for $L(\zeta_p)^+/\Q$ by Theorem \ref{thm:p-Beilinson-abelian}.
Finally, Schneider's conjecture $\Sch(L(\zeta_p), p, r)$ holds for all but finitely
many $r$ by Remark \ref{rmk:Schneider-almost}. Thus the result follows
from Theorem \ref{thm:ETNC-at-r}.
\end{proof}

\begin{remark}
Of course, the result of Corollary \ref{cor:classical-result} is not new. 
In fact, the ETNC for the pair
$(\Q(r)_L, \Z[G])$ holds for every integer $r$ whenever $L/\Q$ is abelian.
If $r \leq 0$ this is the main result of Burns and Greither in \cite{MR1992015}
(important difficulties with the prime $2$ have subsequently been resolved 
by Flach \cite{MR2863902}). The case $r>0$ is due to
Burns and Flach \cite{MR2290586}.
A slightly weaker variant of the ETNC, where the integral group ring $\Z[G]$
is essentially replaced by a maximal order containing it,
has been studied earlier by Huber and Kings \cite{MR2002643}.
\end{remark}

\begin{example} \label{ex:Aff(q)}
Let $E/\Q$ be a Galois extension of totally real fields 
with Galois group $G \simeq \Aff(q)$, where
$q = \ell^n$ is a prime power. Let $r>1$ be an integer.
Since Gross's conjecture holds for all $\chi \in R(G)$ by 
Theorem \ref{thm:Gross-cases} (iv), we have that $T\Omega(\Q(r)_L, \Z[\half][G])$ is rational
by \cite[Theorem 6.5 (i)]{wild-kernels}.
Let us write $\Aff(q) \simeq N \rtimes H$, where $N$ denotes the commutator subgroup
of $\Aff(q)$. Then the $p$-adic group ring $\Z_p[\Aff(q)]$ is `$N$-hybrid'
in the sense of \cite[Definition 2.5]{MR3461042} by \cite[Example 2.16]{MR3461042}
for every prime $p \not= \ell$. Since every $p$-adic group ring is $\{1\}$-hybrid,
we deduce from \cite[Theorem 10.2]{BS-MathZ} 
(or just as well from Theorem \ref{thm:EIMC-known-cases} (ii))
that the equivariant Iwasawa main conjecture
holds unconditionally for $E(\zeta_p)^+_{\infty}/\Q$ for every odd prime $p$.
Thus Conjecture \ref{conj:higher-p-adic-formula} 
holds by Theorem \ref{thm:p-class-holds} 
whenever $\Sch(E(\zeta_p),p,r)$ holds. In particular, this conjecture holds for
almost all $r>1$ for a fixed prime $p$.

We now assume for simplicity that $r \equiv 1 \mod (p-1)$.
The $p$-adic Beilinson conjecture holds for all linear characters of $G$
by Theorem \ref{thm:p-Beilinson-abelian}. As we have already observed in the proof
of Theorem \ref{thm:Gross-cases} (iv) there is only one non-linear character $\chi_{\nl}$
of $\Aff(q)$ which is a $\Z$-linear combination of linear characters and of
$\ind_H^G \mathbbm{1}_H$. Hence it suffices to show Conjecture \ref{conj:p-adic-Beilinson}
for the trivial character $\mathbbm{1}_H$, i.e.~for the trivial extension $E^H/E^H$.
Assuming this we can apply Corollary \ref{cor:ETNC-at-r} to deduce that
the $p$-part of the ETNC for the pair $(\Q(r)_E, e_r\Z[\frac{1}{2}][G])$ holds.
So what is missing here (apart from Schneider's conjecture) is a higher analogue
of Colmez's $p$-adic analytic class number formula \cite{MR922806}
and its complex analytic counterpart. 
(A closer analysis of the proof 
of Theorem \ref{thm:ETNC-at-r} shows that similar observations indeed hold for 
arbitrary $r>1$.) 
\end{example}

\bibliography{p-adic-Beilinson-bib}{}
\bibliographystyle{amsalpha}

\end{document}